\documentclass[12pt]{amsart} 
 
\usepackage{amssymb, epsfig} 
\usepackage{amsmath} 
\usepackage{latexsym} 

\newtheorem{Theorem}{Theorem}[section] 
\newtheorem{Definition}[Theorem]{Definition} 
\newtheorem{Proposition}[Theorem]{Proposition} 
\newtheorem{Lemma}[Theorem]{Lemma}

\newtheorem{Remark}{Remark}[section] 
\newtheorem{Example}{Example}[section] 
 
\newcommand{\R}{\ensuremath{\mathbb{R}}}

\newcommand{\1}{1\hspace{-1.4 mm}1} 
 
\begin{document} 
 
 \title[Existence and uniqueness for dislocation equations]{ 
GLOBAL EXISTENCE RESULTS AND UNIQUENESS FOR DISLOCATION EQUATIONS} 
 
\author[G. Barles, P. Cardaliaguet, O. Ley and R. Monneau] 
{Guy Barles, Pierre Cardaliaguet, Olivier Ley  \\
\& R\'egis Monneau
}

\address{ (G. Barles, O. Ley)  
Laboratoire de Math\'ematiques et Physique Th\'eo\-ri\-que \\
F\'ed\'eration Denis Poisson \\
Universit\'e de Tours \\
Parc de Grandmont, 37200 Tours, France \\ {\tt \{barles,ley\}@lmpt.univ-tours.fr}}
\address{
(P. Cardaliaguet) Universit\'e de Bretagne Occidentale  \\
UFR des Sciences et Techniques, 6 Av. Le Gorgeu \\
BP 809, 29285 Brest, France \\
 {\tt pierre.cardaliaguet@univ-brest.fr}} 
\address{ (R. Monneau) CERMICS, Ecole Nationale des Ponts et Chauss\'ees \\  
6 et 8 avenue Blaise Pascal, Cit\'e Descartes  
Champs-sur-Marne, 77455 Marne-la-Vall\'ee Cedex 2, France \\
 {\tt monneau@cermics.enpc.fr}
}

\begin{abstract} 
We are interested in nonlocal Eikonal Equations arising in the study of the 
dynamics of dislocations lines in crystals. For these nonlocal but also 
non monotone equations, only the existence and uniqueness of Lipschitz and  
local-in-time solutions were available in some particular cases. 
In this paper, we propose a definition of weak solutions for which we are 
able to prove the existence for all time. Then we discuss the uniqueness 
of such solutions in several situations, both in the monotone and 
non monotone case. 
\end{abstract} 
 
 
 
\maketitle 
 
\section{Introduction}

In this article we are interested in the dynamics of defects in crystals,  
called dislocations. The dynamics of these dislocations is the main
microscopic explanation of the macroscopic behaviour of metallic crystals 
(see for instance the physical
monographs Nabarro \cite{N}, Hirth and Lothe \cite{HL}, or Lardner \cite{lardner}
for a mathematical presentation). 
A dislocation is a line moving in a crystallographic plane, called
a slip plane. The typical length of such a dislocation line is of the order
of $10^{-6}m$. Its dynamics is given by a normal velocity proportional to the 
Peach-Koehler force acting on this line.\\

This Peach-Koehler force may have two possible
contributions: the first one is the self-force created by the elastic field generated
by the dislocation line itself (i.e. this self-force is a nonlocal function of the shape of the 
dislocation line); the second one is the force created by
everything exterior to the dislocation line, like the exterior stress
applied on the material, or the force created by other defects. In this
paper, we study a particular model introduced in Rodney et al. \cite{rlf03}.\\

More precisely, if, at time $t$, the dislocation line is the boundary of an open set 
$\Omega_t\subset \R^N$ with $N=2$ for the physical application, the normal 
velocity to the set $\Omega_t$ is given by 
\begin{equation}\label{eq:1} 
V_n= c_0\star \1_{\overline\Omega_t} + c_1    
\end{equation} 
where $\1_{\overline\Omega_t}(x)$ is the indicator function of the set  
$\overline\Omega_t$, which is equal to $1$ if $x\in\overline\Omega_t$ and equal to $0$ 
otherwise. 
The function $c_0(x,t)$ is a kernel  
which only depends on the physical properties 
of the crystal and on the choice of the dislocation 
line whose we follow the evolution.
In the special case of application to dislocations, the kernel $c_0$ 
does not depend on time, but to keep a general setting we allow here a 
dependence on the time variable. 
Here $\star$ denotes the convolution is space, namely 
\begin{equation}\label{def-convolution} 
(c_0(\cdot,t)\star \1_{\overline\Omega_t})(x)= \int_{\R^N} c_0(x-y,t)  
\1_{\overline\Omega_t}(y)dy, 
\end{equation} 
and this term appears to be the Peach-Koehler self-force created by the
dislocation itself, while $c_1(x,t)$ is an additional contribution to the
velocity, created by
everything exterior to the dislocation line. We refer to Alvarez et 
al. \cite{ahlm06} for a detailed  presentation and a derivation of this 
model.\\

We proceed as in the level-set approach to derive an equation for the 
dislocation line. We replace the evolution of a  
set $\Omega_t$ (the strong solution), by the evolution of a function $u$ such  
that $\Omega_t=\left\{u(\cdot,t) >0\right\}$. Roughly speaking the dislocation  
line is represented by the zero level-set of the function $u$ which solves the following equation 
\begin{equation}\label{dislocation} 
\left\{\begin{array}{l} 
\displaystyle{\frac{\partial u}{\partial t}   
= (c_0 (\cdot ,t) \star \1_{\{u(\cdot,t)\geq 0\}}(x) +c_1(x,t))|D u|  
\quad \mbox{in} 
\quad  \R^N\times (0,T)}\\ 
u(\cdot,0)=u_0\quad \mbox{in  }\R^N\; , 
\end{array}\right. 
\end{equation} 
where (\ref{def-convolution}) now reads 
\begin{eqnarray} \label{conv-space} 
c_0 (\cdot ,t) \star \1_{\{u(\cdot,t)\geq 0\}}(x)=\int_{\R^N} c_0(x-y,t)  
\1_{\{u(\cdot,t)\geq 0\}}(y)dy. 
\end{eqnarray} 
Note that (\ref{dislocation}) is not really a level-set equation since it is 
not invariant under nondecreasing changes of functions $u\to \varphi (u)$ where  
$\varphi$ is nondecreasing. As noticed by Slep$\rm \check{c}$ev \cite{slepcev03}, the natural 
level-set equation should be (\ref{slepcev-intro}), see Section \ref{sub-sect-noyau-pos}.\\
 
Although equation (\ref{dislocation}) seems very simple, there are only a few known results. Under suitable assumptions  
on the initial data and on $c_0,c_1$, the existence and uniqueness of the solution is known in two particular cases:  
either for short time (see \cite{ahlm06}), or for all time under the additional assumption that $V_n\ge 0$, which is for instance always satisfied for $c_1$ satisfying $\displaystyle  c_1 (x,t) \geq |c_0(\cdot, 
t)|_{L^1(\R^N)}$ (see \cite{acm05}, \cite{cm06}  or \cite{bl06} for a 
level-set formulation).\\ 
 
In the general case, the existence for all time of solutions to Equation (\ref{dislocation}) is not known and, in particular, in
the case when the kernel $c_0$ has negative values; indeed, in this case, the front propagation problem (\ref{dislocation})
does not satisfy any monotonicity property (preservation of inclusions) and therefore, even if a level-set type equation
can be derived, viscosity solutions' theory cannot be used readily. At this point, it is worth pointing out that a key
property in the level-set approach is the comparison principle for viscosity solutions which is almost equivalent to
this monotonicity property (See for instance Giga's monograph \cite{giga06}). On the other hand, one may try to use
partly viscosity solutions' theory together with some other approximation and/or compactness arguments to prove at least
the existence of weak solutions (in a suitable sense). But here also the bad sign of the kernel creates difficulties since
one cannot use readily the classical half-relaxed limits techniques to pass to the limit in the approximate problems.
Additional arguments are needed to obtain weak solutions.\\

The aim of this paper is to describe a general approach of these dislocations' 
dynamics, based on the level-set approach, which allows us to introduce a 
suitable notion of weak solutions, to prove the existence of these weak  
solutions for all time and to analyse the uniqueness (or non-uniqueness)  
of these solutions. 
 
\subsection{Weak solutions of the dislocation equation}
 
We introduce the following definition of weak solutions, which uses itself 
the definition of $L^1$-viscosity solutions, recalled in Appendix A. 
\begin{Definition}\label{defi:1}{\bf (Classical and weak solutions)}\\ 
For any $T>0$, we say that a function $u\in W^{1,\infty}(\R^N\times [0,T))$ 
is a weak solution of equation (\ref{dislocation}) on the time interval 
$[0,T)$, if there is some measurable map $\chi:\R^N\times (0,T)\to [0,1]$ such that $u$ is a $L^1$-viscosity solution of 
\begin{equation}\label{FormeFaible} 
\left\{\begin{array}{l} 
\displaystyle{\frac{\partial u}{\partial t}  = {\bar c}(x,t) |D u| \quad \mbox{in} 
\quad  \R^N\times (0,T)}\\ 
u(\cdot,0)=u_0\quad \mbox{in  }\R^N\; , 
\end{array}\right. 
\end{equation} 
where  
\begin{equation}\label{eq:chi} 
{\bar c}(x,t) = c_0(\cdot,t)\star \chi(\cdot,t) (x) +c_1(x,t) 
\end{equation} 
and 
\begin{equation}\label{FormeFaible2} 
\begin{array}{l} 
\displaystyle{  \1_{\{u(\cdot,t) > 0\}}(x) \;  \leq \;  \chi(x,t) \; \leq  \1_{\{u(\cdot,t) \geq 0\}}(x) \; ,} 
\end{array} 
\end{equation} 
for almost all $(x,t)\in \R^N\times  [0,T]$. We say that $u$ is a classical solution of equation (\ref{dislocation}) if 
$u$ is a weak solution to (\ref{FormeFaible}) and if 
\begin{equation}\label{FormeFaible3} 
\1_{\{u(\cdot,t) > 0\}}(x) =\1_{\{u(\cdot,t) \geq 0\}}(x) 
\end{equation} 
for almost all $(x,t)\in \R^N\times  [0,T]$. 
\end{Definition} 
\noindent Note that we have 
$\chi(x,t)=\1_{\{u(\cdot,t) > 0\}}(x) =\1_{\{u(\cdot,t) \geq 0\}}(x)$ for 
almost all $(x,t)\in \R^N\times  [0,T]$ 
for classical solutions. 
 
To state our first existence result, we use the following assumptions\\ 
 
\noindent {\bf (H0)} $u_0\in W^{1,\infty}(\R^N)$, $-1 \leq u_0 \leq 1$ and there exists  
$R_0>0$ such that $u_0(x) \equiv -1$  
for $|x|\geq R_0$, \\ 
 
\noindent {\bf (H1)} $c_0\in C([0,T); L^1\left(\R^N \right))$, $D_x c_0\in L^\infty ([0,T); L^1\left(\R^N \right))$,  
$c_1\in C(\R^N\times [0,T))$ and there exists constants $M_1, L_1$ such that, for any $x,y\in \R^N$ and $t\in [0,T]$ 
$$ |c_1(x,t)|\le M_1 \quad \hbox{and} \quad |c_1(x,t)-c_1(y,t)|\le L_1|x-y| .$$ 
 
In the sequel, we denote by $M_0, L_0$, constants such that, for any (or almost every) $t\in [0,T)$, we have 
$$ 
|c_0(\cdot,t)|_{L^1(\R^N)}\le M_0 \quad \hbox{and} \quad 
|D_x c_0(\cdot ,t)|_{L^1(\R^N)} \le L_0 .$$ 
 
\smallskip 
Our first main result is the following. 
\begin{Theorem}\label{th:1}{\bf (Existence of weak solutions)}\\ 
Under assumptions {\bf (H0)}-{\bf (H1)}, for any $T >0$ and for any initial data $u_0$,  
there exists a weak solution of equation (\ref{dislocation}) on the time 
interval $[0,T)$ in the sense of Definition \ref{defi:1}. 
\end{Theorem} 
 
\smallskip 
 
Our second main result states that a weak solution is a classical one if the evolving set is expanding and if the 
following additional condition is fulfilled\\ 
 
\noindent {\bf (H2)} $c_1$ and $c_0$ satisfy {\bf (H1)} and there exists constants $m_0, N_1$ and a positive function $N_0\in L^1(\R^N)$ such that, for any $ x,h\in \R^N$, $t\in [0,T)$, we have 
$$|c_0(x,t)| \le m_0 ,$$ 
$$ |c_1(x+h,t)+c_1(x-h,t)-2c_1(x,t)|\le N_1 |h|^2,  $$ 
$$ 
|c_0(x+h,t)+ c_0(x-h,t)-2c_0(x,t)| \le N_0(x) |h|^2 .$$ 
 
\smallskip 
 
\begin{Theorem}\label{th:2}{\bf (Some links between weak solutions and 
classical continuous viscosity solutions and uniqueness results)}\\ 
Assume {\bf (H0)}-{\bf (H1)} and suppose that there is some $\delta\geq 0$ such that,
for all measurable map $\chi:\R^N\times (0,T)\to [0,1],$ 
\begin{equation}\label{HypPos} 
for \ all \ (x,t)\in \R^N \times [0,T],\quad  
c_0 (\cdot ,t) \star \chi (\cdot ,t) (x) +c_1(x,t)\geq \delta, 
\end{equation} 
and that the initial data $u_0$ satisfies (in the viscosity sense) 
\begin{eqnarray} \label{lgb1111} 
- |u_0| - |D u_0| \leq -\eta_0 \quad\hbox{in  } \R^N \; , 
\end{eqnarray} 
for some $\eta_0 >0$. Then any weak solution $u$ of 
(\ref{dislocation}) in the sense of Definition \ref{defi:1}, is a classical 
continuous viscosity solution of (\ref{dislocation}).  
This solution is unique if  {\bf (H2)} holds and
 
(i) either $\delta>0$, 
 
(ii) or $\delta=0$ and $u_0$ is semiconvex, i.e. satisfies for some constant $C>0$: 
$$u_0(x+h)+u_0(x-h)-2u_0(x)  \ge -C |h|^2, \quad \forall x,h\in \R^N.$$ 
\end{Theorem} 
Assumption (\ref{HypPos}) ensures that the velocity $V_n$ in (\ref{eq:1}) is 
positive for positive $\delta$. Of course, we can state similar results in the case of negative velocity.
Assumption (\ref{lgb1111}) means that $u_0$ is a viscosity 
subsolution of $-|v(x)|-|D v(x)|+\eta_0\leq 0.$  When  $u_0$ is $C^1,$ it 
follows that the gradient of $u_0$ does not vanish on the set $\{u_0=0\}$  
(see  \cite{ley01} for details). 
Point (ii) of the theorem is the main result of \cite{acm05, bl06}. 
We also point out that, with adapted proofs, only a bound from below could be required in  
{\bf (H2)} on $c_1(x+h,t)+c_1(x-h,t)-2c_1(x,t)$ and  $c_0(x+h,t)+ c_0(x-h,t)-2c_0(x,t)$. 

\begin{Remark} \label{unicite-c0pos1}  \rm 
In particular
Theorem \ref{th:2} implies uniqueness in the case
$c_0\geq 0$ and $c_1\equiv 0.$ 
The general study of nonnegative kernels is provided below
(see Theorem \ref{th:4} and Remark \ref{rem-sl1}).
\end{Remark} 
 
\subsection{Nonnegative kernel $c_0\ge 0$} 
\label{sub-sect-noyau-pos} 
 
In the special case where the kernel $c_0$ is non-negative, an inclusion 
principle for the dislocations lines, or equivalently a comparison principle 
for the  functions of the level-set formulations is expected (cf. Cardaliaguet \cite{cardaliaguet00} and Slep$\rm \check{c}$ev 
\cite{slepcev03}).  
 
Moreover, in the classical level-set approach, all the level-sets of $u$ should have the  
same type of normal velocity and Slep$\rm \check{c}$ev 
\cite{slepcev03} remarked that a formulation with a nonlocal term of the form  
$\{u(\cdot,t)\geq u(x,t)\}$ is more appropriate.  
Therefore it is natural to start studying the following equation (which replaces 
equation (\ref{dislocation})) 
\begin{equation}\label{slepcev-intro} 
\left\{\begin{array}{l} 
\displaystyle{\frac{\partial u}{\partial t}   
= (c_0 (\cdot ,t) \star \1_{\{u(\cdot,t)\geq u(x,t)\}}(x) +c_1(x,t))|D u|  
\quad \mbox{in} 
\quad  \R^N\times (0,T) \; ,}\\ 
u(\cdot,0)=u_0\quad \mbox{in  }\R^N\; , 
\end{array}\right. 
\end{equation} 
where $\star$ denotes the convolution in space as in (\ref{conv-space}). 
 
The precise meaning of a viscosity solution of (\ref{slepcev-intro}) is 
given in Definition \ref{Defi-slepcev}.  
 
In this context, assumption {\bf (H0)} can be weakened into the following condition which allows to consider unbounded evolving sets\\ 
 
\noindent {\bf (H0')} $u_0\in BUC(\R^N)$.\\ 
 
Our main result for this 
equation is 
\begin{Theorem}\label{th:3}{\bf (Existence and uniqueness)}\\ 
Assume that $c_0\ge 0$ on $\R^N\times [0,T]$ and that 
{\bf (H0')}-{\bf (H1)} hold. Then there exists a unique viscosity solution $u$ of (\ref{slepcev-intro}). 
\end{Theorem} 
 
\begin{Remark} \label{rem-slep0} \rm The comparison principle for this equation (see Theorem \ref{comp-slep}) is a generalization of 
\cite[Theorem 2.3]{slepcev03}: indeed, in \cite{slepcev03}, everything takes 
place in a fixed bounded set whereas here one has to deal with unbounded 
sets. See also \cite{dlfm07} for related results. 
\end{Remark} 
 
Now we turn to the connections with weak solutions. To do so, if $u$ is the unique continuous solution of (\ref{slepcev-intro}) given by Theorem \ref{th:3},  
we introduce the functions $\rho^+, \rho^-: \R^N \times [0, T] \to \R$ defined by 
$$ \rho^+ : =\1_{\{u\geq 0\}}\quad \hbox{and} \quad \rho^- : =\1_{\{u > 0\}}\; .$$ 
 
Our result is the 
\begin{Theorem}\label{th:4}{\bf (Maximal and minimal weak solutions)}\\ 
Under the assumptions of Theorem~\ref{th:3}, the maximal and minimal weak solutions of (\ref{dislocation}) are the continuous functions $v^+, v^-$ which are the unique $L^1$-viscosity solutions of the equations 
\begin{equation}\label{nonlocaldisc} 
\left\{\begin{array}{l} 
\displaystyle{ \frac{\partial v^\pm}{\partial t} =  c[\rho^\pm](x,t)|D v^\pm| \quad \mbox{in} 
\quad  \R^N\times (0,T)\; ,} \\[2mm] 
v^\pm (x,0)=u_0 (x) \quad \mbox{in} 
\quad  \R^N\; , 
\end{array} 
\right. 
\end{equation} 
where 
$$c[\rho](x,t) := c_0 (\cdot ,t) \star \rho(\cdot,t) (x) +c_1(x,t)\quad \mbox{in} 
\quad  \R^N\times (0,T)\; .$$ 
The functions $v^\pm$ satisfy $\{v^+ (\cdot, t) \geq 0\} = \{u(\cdot, t)\geq 0\}$ and  
$\{v^- (\cdot, t) > 0\} = \{u(\cdot, t)> 0\}$, where $u$ is the solution  
of (\ref{slepcev-intro}). 

Moreover, if the set $\{u(\cdot, t) = 0\}$  
has a zero-Lebesgue measure for almost all $t \in (0,T)$, then problem (\ref{dislocation}) 
has a  unique weak solution which is also a classical one.   
\end{Theorem} 
 
\begin{Remark} \label{rem-sl1} {\rm  \ \\ 
1.  Theorem~\ref{th:4} shows that, in the case when $c_0 \geq 0$, Slep$\rm\check{c}$ev's approach  
allows to identify the maximal and minimal weak solutions as being associated 
to $\rho^\pm$. \\ 
\noindent 2. Equalities $\{v^- (\cdot, t) \geq 0\} = \{u(\cdot, t)\geq 0\}$ 
and $\{v^+ (\cdot, t) > 0\} = \{u(\cdot, t)> 0\}$ do not hold in general  
(see for instance example \ref{exemple00} in Section \ref{ws}).\\ 
\noindent 3. If the set $\{v^\pm (\cdot, t) = 0\}$ develops an interior, a dramatic loss 
of uniqueness for the weak solution of (\ref{dislocation}) may occur.  
This is illustrated by Example \ref{exemple00} below, where we are able to built infinitely many solutions 
after the onset of fattening. \\ 
\noindent 4. We have uniqueness for (\ref{dislocation}) if $\{u(\cdot, t) = 0\}$  
has a zero-Lebesgue measure for almost all $t \in (0,T).$ This condition is fulfilled
when, for instance, $c[\rho]\geq 0$ holds for any 
indicator function $\rho$ and (\ref{lgb1111}) holds (see also Remark
\ref{unicite-c0pos1}).  
}
\end{Remark} 
 
\subsection{Organization of the paper} 
 
In Section~\ref{bree}, we recall basic results for the classical Eikonal 
Equation which are used throughout the paper.
In Section~\ref{ws}, we prove the existence of weak solutions for 
equation (\ref{dislocation}), namely Theorem~\ref{th:1} and give a 
counter-example to the uniqueness in general. 
Let us mention that this first part of the paper, even if it requires rather deep results of viscosity solutions theory, is of a general interest for a wide audience and can be read without having an expertise in this theory since one just need to apply the results which are anyway rather natural.
In Section~\ref{sect4}, we prove Theorem~\ref{th:2} in the case of expanding dislocations. 
The arguments we use here are far more involved from a technical point of view: in particular we need some fine estimates of the perimeter of the evolving sets. In Section 
~\ref{slep}, we study the Slep$\rm \check{c}$ev formulation in the case of 
non-negative kernels, and prove Theorems \ref{th:3} and \ref{th:4}. 
In the spirit, this section is closely related to the classical level-set
approach but is more technical.
Finally, 
for sake of completeness, we recall in Appendix A the Definition of $L^1$-viscosity 
solutions and a new stability result proved by Barles in \cite{barles06}.

\section{Some basic results for the classical (local) eikonal equation}\label{bree} 
 
We want to recall in this section some basic results on the level-set equation 
\begin{equation}\label{eq:2} 
\left\{\begin{array}{l} 
\displaystyle{\frac{\partial v}{\partial t}  = a (x,t)|D v| \quad \mbox{in} 
\ \R^N\times (0,T)}\\ 
\\ 
v(\cdot,0)=u_0\quad \mbox{in} \  \R^N\; , 
\end{array}\right. 
\end{equation} 
where $T >0$ and $a : \R^N\times [0,T] \to \R$ is, at least, a continuous function.  
 
We provide some classical estimates on the solutions to (\ref{eq:2}) when $a$ satisfies suitable assumptions. 
Our result is the following. 
\begin{Theorem}\label{basicthm} If $u_0$ satisfies {\bf (H0)} and $a$ satisfies the assumptions of $c_1$ in {\bf (H1)},  
then Equation (\ref{eq:2}) 
has a unique continuous solution $v$ which is Lipschitz continuous in  
$\R^N\times [0,T]$  
and which satisfies 
\begin{itemize} 
\item[$(i)$] $-1 \leq v \leq 1$ in $\R^N\times (0,T)$, $v (x,t) \equiv -1$ for  
$|x|\geq R_0 + M_1t$,  
\item[$(ii)$] $ | Dv(\cdot ,t)|_\infty \leq | Du_0|_\infty \, {\rm e}^{L_1t}$, 
\item[$(iii)$] $\displaystyle{ \left|v_t(\cdot ,t)\right|_\infty \leq M_1| Du_0|_\infty \, {\rm e}^{L_1t}\, } 
.$ 
\end{itemize} 
\end{Theorem} 
 
We skip the very classical proof of Theorem~\ref{basicthm}; we just point out  
that the first point comes from the  
comparison result for (\ref{eq:2}) and the ``finite speed of propagation  
property'' (See Crandall \& Lions \cite{cl83}),  
while the second one is a basic gradient estimate (see for example Ley  
\cite{ley01}) and the last one comes directly  
from the fact that the equation is satisfied almost everywhere. 
 
The main consequence of this result is that the solution remains in a compact  
subset of the Banach space 
$(C(\R^N\times [0,T]), |\cdot |_\infty)$ as long  
as $u_0$ and $a$ satisfies {\bf (H0)}-{\bf (H1)} with fixed constants.\\ 
 
Let us introduce the following 
\begin{Definition}{\bf (Interior ball property)}\\ 
We say that a closed set $K \subset \R^N$ has an interior ball property 
of radius $r>0$, if for any $x\in K$, there exists $p\in 
\R^N\backslash \left\{0\right\}$ such that $B(x-r \frac{p}{|p|},r) \subset K$. 
\end{Definition} 
 
We will also use the following result, due to Cannarsa and Frankowska \cite{cf06}, the proof of which is given in Appendix B for sake of completeness. 
\begin{Lemma}\label{lem:bornebouleint}{\bf (Interior ball 
    regularization)}\\ 
Suppose  {\bf (H0)} and that $a$ satisfies the assumptions of $c_1$ in {\bf (H1)}-{\bf (H2)} and  
there exists a constant $\delta>0$ such that  
$$c_1 \ge \delta >0\quad\mbox{on}\quad \R^N\times [0,T].$$ 
Then there exists a constant $\gamma$ (depending in particular on $\delta >0$ 
and $T$ and on the other constants of the problem) such that for the solution 
$v$ of (\ref{eq:2}), the set $\left\{v(\cdot,t)\ge 0\right\}$ has an 
interior ball property of radius $r_t \ge \gamma t$ for $t\in (0,T)$. 
\end{Lemma} 

\section{Existence of weak solutions for equation (\ref{dislocation})}\label{ws} 
We aim at solving equation (\ref{dislocation}), i.e. 
$$\left\{\begin{array}{l} 
\displaystyle{\frac{\partial u}{\partial t}   
= (c_0 (\cdot ,t) \star \1_{\{u(\cdot,t)\geq 0\}}(x) +c_1(x,t))|D u|  
\quad \mbox{in} 
\quad  \R^N\times (0,T)}\\ 
u(\cdot,0)=u_0\quad \mbox{in  }\R^N\; , 
\end{array}\right.$$ 
proving Theorem \ref{th:1}, which states the existence of weak solutions  
as introduced in Definition \ref{defi:1}.

A key difficulty to solve (\ref{dislocation}) comes from the fact that, in 
this kind of level-set equations, one may face  
the so-called ``non-empty interior difficulty'', i.e. that the $0$--level-set of  
the solution is ``fat'' which  
may mean either that it has a non-empty interior or a non-zero Lebesgue measure.  
Clearly, in both cases,  
$\1_{\{u(\cdot,t)\geq 0\}}$ is different from $\1_{\{u(\cdot,t) > 0\}}$ and this  
leads to rather bad stability  
properties for (\ref{dislocation}) and therefore to difficulties to prove the existence  
of a solution (and even more for the  
uniqueness). The notion of weak solution (\ref{FormeFaible})-(\ref{eq:chi})-(\ref{FormeFaible2}) emphasizes  
this difficulty. On the contrary, if ${\bar c}(x,t) \geq 0$ in  
$\R^N \times [0,T]$, it is known  
that the ``non-empty interior difficulty'' cannot happen (See Barles, Soner and 
Souganidis \cite{bss93} and  
Ley \cite{ley01}) and we recover a more classical formulation. We discuss this question in the next section, as well as some uniqueness 
issues for our weak solutions. Let us finally note that weak solutions for 
(\ref{dislocation}) satisfy the following inequalities: 
\begin{Proposition}\label{Trestresweak} Let $u$ be a weak solution to (\ref{dislocation}).  
Then $u$ also satisfies in the $L^1$-sense 
\begin{equation}\label{sub} 
\frac{\partial u}{\partial t}  \leq \left( c^+_0 (\cdot,t)\star 
  \1_{\{u(\cdot,t)  \geq  
0\}} (x) 
- c^-_0 (\cdot,t) \star \1_{\{u(\cdot,t) > 0\}}(x) +c_1(x,t) \right) |D u|, 
\end{equation} 
\begin{equation}\label{super} 
\frac{\partial u}{\partial t}  \geq \left( c^+_0 (\cdot,t)\star \1_{\{u(\cdot,t) > 0\}}(x)  
- c^-_0(\cdot,t) \star \1_{\{u(\cdot,t) \geq 0\}}(x) +c_1(x,t) \right) |D u|,  
\end{equation} 
in $\R^N\times (0,T),$ 
where $c_0^+=\max(0,c_0)$ and $c_0^-=\max(0,-c_0)$. 
\end{Proposition} 
 
\noindent{\bf Proof of Proposition \ref{Trestresweak}.} \\ 
Let  ${\bar c}$ be associated with $u$ as in (\ref{FormeFaible})-(\ref{eq:chi})-(\ref{FormeFaible2}). Then we have 
$$ 
{\bar c}(x,t)  \geq  c^+_0 (\cdot,t)\star \1_{\{u(\cdot,t) > 0\}}(x)  
- c^-_0 (\cdot,t)\star \1_{\{u(\cdot,t) \geq 0\}}(x) +c_1(x,t)  
$$ 
for every $x\in \R^N$ and almost every $t\in (0,T)$. We note that the right-hand side of the inequality is lower-semicontinuous. 
Following Lions and Perthame \cite{lp87}, $u$ then solves (\ref{super}) in the usual viscosity sense. The proof of (\ref{sub}) can be achieved 
in a similar way. 
\hfill $\Box$ \\

\noindent{\bf Proof of Theorem \ref{th:1}.} \\ 
\noindent{\it 1. Introduction of a perturbated equation.}  
First we are going to solve the equation 
\begin{eqnarray}\label{pert-eq} 
\hspace*{0.8cm} 
\frac{\partial u}{\partial t}  = \left( c_0 (\cdot,t)\star \psi_\varepsilon 
  (u (\cdot,t))  (x)+c_1(x,t)  
\right) |D u|  
\quad \mbox{in  }  \R^N\times (0,T)\; ,  
\end{eqnarray} 
where $\psi_\varepsilon :\R\to \R$ is a sequence of continuous functions such  
that  
$\psi_\varepsilon (t) \equiv 0$ for $t \leq -\varepsilon$,  $\psi_\varepsilon  
(t) \equiv 1$ for $t \geq 0$  
and $\psi_\varepsilon$ is an affine function on $[-\varepsilon, 0]$. 
 
We aim at applying Schauder's fixed point Theorem to a suitable map. We note that an alternative proof could be given by using techniques 
developed by Alibaud in \cite{alibaud07}. \\ 
 
\noindent{\it 2. Definition of a map $\mathcal{T}.$}  
We introduce the  convex and compact (by Ascoli's Theorem) subset 
\begin{eqnarray*} 
X= \{ u\in C(\R^N\times [0,T]) : u\equiv -1 \ {\rm in} \ \R^N \backslash  
B(0,R_0 +MT), \\ 
\hspace*{4cm} |Du|, |u_t|/M \le |Du_0|_\infty  e^{L T}\} 
\end{eqnarray*} 
of $(C(\R^N\times [0,T]), |\cdot |_\infty)$ for $M=M_0+M_1$ and $L=L_0+L_1$,  
and the map ${\mathcal T} : X\to X$ defined by~:  
if $u\in C(\R^N\times [0,T])$, then ${\mathcal T} (u)$ is the unique solution $v$ of  
(\ref{eq:2}) for  
\begin{eqnarray*} 
c_\varepsilon (x,t) &=& c_0 (\cdot,t)\star \psi_\varepsilon (u (\cdot,t)) (x) +c_1(x,t) \\ 
&=& \int_{\R^N} c_0(x-z,t)  \psi_\varepsilon (u (z,t))dz+c_1(x,t). 
\end{eqnarray*} 
This definition is justified by the fact that, under assumption {\bf (H1)} 
on $c_1$ and $c_0$,  
$c_\varepsilon$ satisfies {\bf (H1)} with fixed constants $M=M_0+M_1$ and $L=L_0+L_1$;  
indeed $M$ is a bound on $\sup_{[0,T]}\, |c_0(\cdot, t)|_{L^1} + M_1$ while $L$  
is estimated by the following calculation: for all $x,y\in \R^N,$ $t\in [0,T]$  
and $u\in X$, we have 
\begin{eqnarray} \label{calcul-equicont} 
&& c_\varepsilon (x,t)- c_\varepsilon (y,t) \\ 
&=& \int_{\R^N} (c_0(x-z,t)-c_0(y-z,t)) \psi_\varepsilon (u(z,t))dz +c_1(x,t)- 
c_1(y,t)  \nonumber\\ 
&\leq & 
\int_{\R^N} |c_0(x-z,t)-c_0(y-z,t)|dz + |c_1(x,t)-c_1(y,t)| \nonumber \\ 
&\leq & (L_0 + L_1)|x-y|, \nonumber 
\end{eqnarray} 
since $0\leq \psi_\varepsilon\leq 1$. 
 
Finally, under assumption {\bf (H0)}-{\bf (H1)}, for any $u\in X,$ the results of Theorem \ref{basicthm} apply to  
(\ref{pert-eq}) which imply that ${\mathcal T} (u) \in X$. 
It follows that $\mathcal{T}$ is well-defined. \\ 
 
\noindent{\it 3. Application of Schauder's fixed point Theorem to  
$\mathcal{T}.$}  
The map ${\mathcal T}$ is continuous since $\psi_\varepsilon$ is continuous, by  
using the classical  
stability result for viscosity solutions (see for  
instance (\ref{eq:stab0}) in Section \ref{sect4}).  
Therefore ${\mathcal T}$ has a fixed point $u_\varepsilon$ which is bounded in  
$W^{1,\infty} (\R^N\times [0,T])$  
uniformly with respect to $\varepsilon$ (since $M$ and $L$ are independent of  
$\varepsilon$). \\ 
 
\noindent{\it 4. Convergence of the fixed point when $\varepsilon\to 0.$}  
From Ascoli's Theorem, we extract a subsequence  
$(u_{\varepsilon'})_{\varepsilon'}$ which converges locally uniformly  
to a function denoted by $u$ (in fact globally since the $u_{\varepsilon'}$ are  
equal to $-1$ outside a fixed compact subset).  
 
The functions $\chi_{\varepsilon'}:=\psi_{\varepsilon'} (u_{\varepsilon'})$ satisfy $0\le\chi_{\varepsilon'} \le 1$. Therefore we can extract a 
subsequence---still denoted $(\chi_{\varepsilon'})$---which  converges 
weakly$-*$ in $L^\infty_{\rm loc}(\R^N\times [0,T])$ to some function $\chi:\R^N\times (0,T)\to [0,1]$. 
Therefore, for all $\varphi\in L^1_{\rm loc}(\R^N\times [0,T]),$ 
\begin{eqnarray} \label{cvLinf} 
\int_0^T \int_{\R^N} \varphi \chi_{\varepsilon'} dxdt \to \int_0^T \int_{\R^N} \varphi \chi dx dt. 
\end{eqnarray} 
From Fatou's lemma, if $\varphi$ is nonnegative, it follows 
\begin{eqnarray*} 
\int_0^T \int_{\R^N} \varphi  (x,t)\chi  (x,t)dxdt &\leq & \int_0^T \int_{\R^N} \varphi (x,t) 
\mathop{\rm lim\,sup}_{\varepsilon'\to 0} \chi_{\varepsilon'} (x,t) dxdt \\ 
&\leq & 
\int_0^T \int_{\R^N} \varphi (x,t) 
\mathop{\rm lim\,sup}_{\varepsilon'\to 0, x'\to x, t'\to t} 
\chi_{\varepsilon'} (x',t') dxdt \\ 
&\leq & 
\int_0^T \int_{\R^N} \varphi (x,t) \1_{\{u(\cdot,t) \geq 0\}}(x)dxdt. 
\end{eqnarray*} 
Since the previous inequalities hold for any nonnegative $\varphi\in L^1_{\rm loc}(\R^N\times [0,T]),$ 
we obtain that, for  almost every $(x,t)\in \R^N\times (0,T),$ 
$$ 
\chi(x,t) \leq \1_{\{u(\cdot,t) \geq 0\}}(x). 
$$ 
Similarly we get 
$$ 
\displaystyle{\1_{\{u(\cdot,t) > 0\}}(x)} \le \chi(x,t). 
$$ 
Furthermore, setting $c_{\varepsilon'}= c_0 \star \chi_{\varepsilon'} + c_1,$ 
from (\ref{cvLinf}), we have, for all $(x,t)\in \R^N\times [0,T],$ 
\begin{eqnarray*} 
\int_0^t c_{\varepsilon'}(x,s)ds  
&= & \int_0^t \int_{\R^N} c_0(x-y,s)\chi_{\varepsilon'} (y,s)dyds + \int_0^t c_1(x,s)ds \\ 
& \to & \int_0^t {\bar c}(x,s)ds, 
\end{eqnarray*} 
where ${\bar c}(x,t)=c_0(\cdot, t)\star \chi(\cdot,t) (x) + c_1(x,t).$ 
The above convergence is pointwise but, noticing that $c_{\varepsilon'}$ 
satisfies ${\bf (H3)}$ (with $M:=M_0+M_1$ and $L:=L_0+L_1$) and using  
Remark \ref{rem-cvunif}, we can apply the stability Theorem \ref{L1stability} 
given in Appendix A.  
We obtain that  
$u$ is $L^1$-viscosity solution to (\ref{FormeFaible}) with ${\bar c}$ satisfying (\ref{eq:chi})-(\ref{FormeFaible2}). 
\hfill $\Box$ \\ 
 
The following example is inspired from \cite{bss93}.  
 
\begin{Example} \rm (Counter-example to the uniqueness of weak solutions)\label{exemple00}\\ 
Let us consider, in dimension $N=1$, the following equation of type 
(\ref{dislocation}), 
\begin{equation}\label{dislocation-d1} 
\left\{\begin{array}{l} 
\displaystyle{\frac{\partial U}{\partial t}   
= (1 \star \1_{\{U(\cdot,t)\geq 0\}}(x) +c_1(t))|D U|  
\quad \mbox{in} 
\quad  \R\times (0,2]}\\[2mm] 
U(\cdot,0)=u_0\quad \mbox{in  }\R\; , 
\end{array}\right. 
\end{equation} 
where we set 
$c_0(x,t):= 1,$ $c_1(x,t):=c_1(t)= 2(t-1)(2-t)$ and $u_0(x)=1-|x|.$ 
Note that $1 \star \1_A ={\mathcal L}^1 (A)$ for any measurable 
set $A\subset \R,$ where ${\mathcal L}^1 (A)$ is the Lebesgue measure 
on $\R.$ 
 
We start by solving auxiliary problems for time in $[0,1]$ and $[1,2]$ in order to 
produce a family of solutions for the original problem in $[0,2].$ 
 
\noindent 
{\it 1. Construction of a solution for $0\leq t\leq 1.$} The function $x_1(t)=(t-1)^2$ is the solution of 
the ode 
\begin{eqnarray*} 
\dot{x}_1 (t)= c_1(t)+2x_1(t) \ {\rm for} \ 0\leq t\leq 1, \ \ \ 
{\rm and } \ x(0)=1, 
\end{eqnarray*}  
(note that $\dot{x}_1\leq 0$ in $[0,1]$). 
Consider 
\begin{eqnarray} \label{eq-eik-1} 
\left\{\begin{array}{l} 
\displaystyle{\frac{\partial u}{\partial t}= \dot{x}_1 (t)\left|\frac{\partial 
  u}{\partial x}\right| \quad\hbox{in  } \R \times (0,1],} \\ 
u(\cdot,0)=u_0\quad\hbox{in  } \R. 
\end{array}\right. 
\end{eqnarray}  
There exists a unique continuous viscosity solution $u$ of (\ref{eq-eik-1}). 
Looking for $u$ under the form $u(x,t)=v(x,\Gamma (t))$ with $\Gamma (0)=0,$ 
we obtain that $v$ satisfies 
\begin{eqnarray*} 
\frac{\partial v}{\partial t}\, \dot{\Gamma}(t)= \dot{x}_1 (t)\left|\frac{\partial 
  v}{\partial x}\right|. 
\end{eqnarray*}  
Choosing $\Gamma (t)=-x_1(t)+1,$ we get that $v$ is the solution of 
\begin{eqnarray*}  
\left\{\begin{array}{l} 
\displaystyle{\frac{\partial v}{\partial t}= -\left|\frac{\partial 
  v}{\partial x}\right| \quad\hbox{in  } \R \times (0,1],} \\ 
v(\cdot,0)=u_0\quad\hbox{in  } \R. 
\end{array}\right. 
\end{eqnarray*}  
By Oleinik-Lax formula, $\displaystyle v(x,t)=\mathop{\rm inf}_{|x-y|\leq t} 
u_0(y).$ Since $u_0$ is even, we have, for all $(x,t)\in\R\times [0,1],$ 
\begin{eqnarray*} 
u(x,t)= \mathop{\rm inf}_{|x-y|\leq \Gamma (t)} u_0(y) =u_0 (|x|+\Gamma (t))= 
u_0 (|x|-x_1(t)+1). 
\end{eqnarray*}  
Therefore, for $0\leq t\leq 1,$ 
\begin{equation} \label{incl678} 
\{ u(\cdot ,t)>0\}= (-x_1(t),x_1(t)) \  \ \ {\rm and} \  \ \ 
\{ u(\cdot ,t)\geq 0\}= [-x_1(t),x_1(t)]. 
\end{equation}  
We will see in Step 3 that $u$ is a solution of (\ref{dislocation-d1}) in $[0,1].$ 
 
\noindent 
{\it 2. Construction of solutions for $1\leq t\leq 2.$} Consider now, for any measurable function $0\leq \gamma (t)  \leq 1,$ the unique 
solution $y_\gamma$ of the ode 
\begin{eqnarray}\label{edo123} 
\dot{y}_\gamma (t)= c_1(t)+2\gamma (t)  y_\gamma (t) \ {\rm for} \ 1\leq t\leq 2, \ \ \ 
{\rm and } \  y_\gamma (1)=0. 
\end{eqnarray}  
By comparison, we have $0\leq y_0(t) \leq y_\gamma (t)\leq y_1 (t)$ for $1\leq 
t\leq 2,$ where $y_0, y_1$ are the solutions of (\ref{edo123}) obtained with $\gamma (t)\equiv 0,1.$ 
In particular, it follows that $\dot{y}_\gamma \geq 0$ in $[1,2].$ 
Consider 
\begin{eqnarray*}  
\left\{\begin{array}{l} 
\displaystyle{\frac{\partial u_\gamma}{\partial t}= \dot{y}_\gamma (t)\left|\frac{\partial 
  u_\gamma}{\partial x}\right| \quad\hbox{in  } \R \times (1,2],} \\ 
u_\gamma(\cdot,1)=u (\cdot ,1)\quad\hbox{in  } \R, 
\end{array}\right. 
\end{eqnarray*}  
where $u$ is the solution of (\ref{eq-eik-1}). 
Again, this problem has a unique continuous viscosity solution $u_\gamma$ and 
setting 
$\Gamma_\gamma (t)=y_\gamma (t)\geq 0$ for $1\leq t\leq 2,$ we obtain that 
$v_\gamma$ defined by $v_\gamma (x,\Gamma_\gamma (t))=u_\gamma (x,t)$ is the 
unique continuous viscosity solution of 
\begin{eqnarray*}  
\left\{\begin{array}{l} 
\displaystyle{\frac{\partial v_\gamma}{\partial t}= \left|\frac{\partial 
  v_\gamma}{\partial x}\right| \quad\hbox{in  } \R \times (0,\Gamma_\gamma (2)],} \\ 
v_\gamma(\cdot,0)=u (\cdot ,1)\quad\hbox{in  } \R. 
\end{array}\right. 
\end{eqnarray*}  
Therefore, for all $(x,t)\in\R\times [1,2],$ we have 
\begin{eqnarray*} 
u_\gamma(x,t)= \mathop{\rm sup}_{|x-y|\leq y_\gamma (t)} u(y,1) = 
\left\{\begin{array}{ll} 
0 & {\rm if} \ |x|\leq y_\gamma (t), \\ 
u(|x|-y_\gamma (t), 1) & {\rm otherwise}. 
\end{array}\right. 
\end{eqnarray*}  
(Note that $u(-x,t)=u(x,t)$ since $u_0$ is even and, since $u(\cdot ,1)\leq 
0,$ by the maximum principle, we have $u_\gamma \leq 0$ in $\R\times [1,2].$) 
It follows that, for all $1\leq t\leq 2,$ 
\begin{equation} \label{incl679} 
\{ u_\gamma(\cdot ,t)>0\}= \emptyset \  \ \ {\rm and} \  \ \ 
\{ u_\gamma(\cdot ,t)\geq 0\}= \{ u_\gamma (\cdot ,t)=0\} = [-y_\gamma (t),y_\gamma (t)]. 
\end{equation}  
 
\noindent 
{\it 3. There are several weak solutions of (\ref{dislocation-d1}).} 
Set, for $0\leq \gamma (t)  \leq 1,$ 
\begin{eqnarray*} 
\begin{array}{lll} 
c_\gamma (t)= c_1(t)+2 x_1(t), & U_\gamma (x,t)= u(x,t) &  {\rm if} \ (x,t)\in \R\times [0,1], \\  
c_\gamma (t)= c_1(t)+2 \gamma (t) y_\gamma (t), & U_\gamma (x,t)= u_\gamma (x,t) &  {\rm if} \ (x,t)\in \R\times [1,2]. \\ 
\end{array} 
\end{eqnarray*}  
Then, from Steps 1 and 2,  $U_\gamma$ is the unique continuous viscosity 
solution of  
\begin{eqnarray} \label{eik567}  
\left\{\begin{array}{l} 
\displaystyle{\frac{\partial U_\gamma}{\partial t}=  c_\gamma (t) \left|\frac{\partial 
  U_\gamma}{\partial x}\right| \quad\hbox{in  } \R \times (0,2],} \\ 
U_\gamma(\cdot,0)=u_0 \quad\hbox{in  } \R. 
\end{array}\right. 
\end{eqnarray}  
Taking $\chi_\gamma (\cdot ,t)=\gamma (t) \1_{[ -y_{\gamma}(t), y_{\gamma} (t) ]}$ for $1\leq t\leq 2,$ 
from (\ref{incl678}) and (\ref{incl679}), we have 
\begin{eqnarray*} 
\1_{\{U_\gamma (\cdot ,t)>0\}}\leq \chi_\gamma (\cdot ,t) \leq \1_{\{U_\gamma (\cdot ,t)\geq 0\}}, 
\end{eqnarray*}  
(see Figure \ref{dess-gamma}). 
It follows that all the $U_\gamma$'s, for measurable $0\leq \gamma (t) \leq 1,$ are all weak 
solutions of (\ref{dislocation-d1}) so we do not have uniqueness and the set of solutions is quite large.  
\begin{figure}[ht] 
\begin{center} 
\epsfig{file=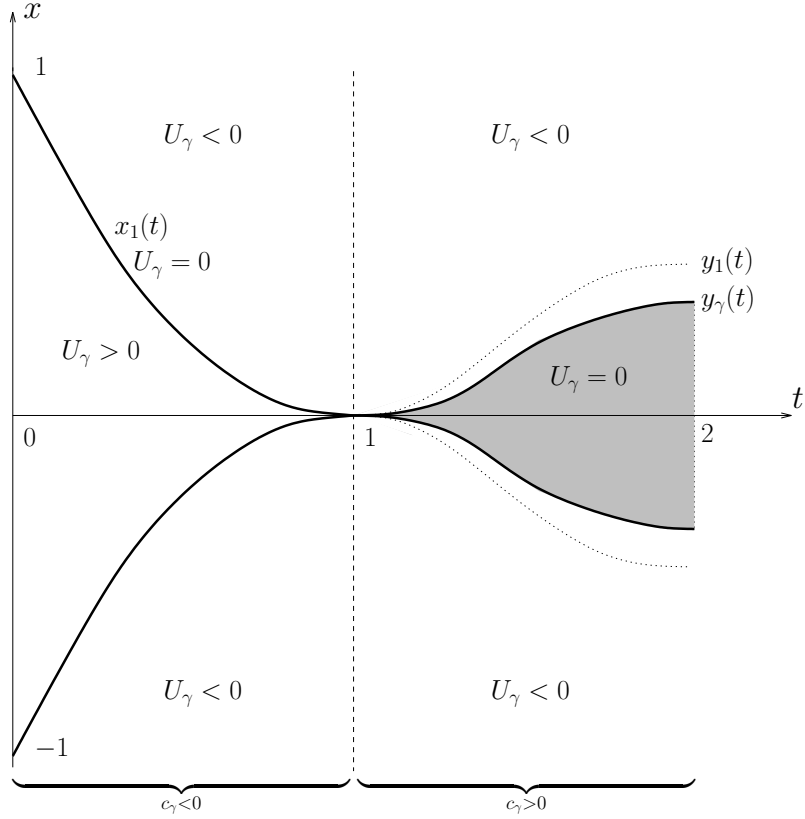, width=10cm}  
\end{center} 
\caption{\label{dess-gamma} 
{\it Fattening phenomenon for the functions $U_\gamma.$}} 
\end{figure} 
 
Let us complete this counter-example by  pointing out: \\ 
(i) as in \cite{bss93}, non uniqueness comes from the fattening phenomenon for 
the front which is due to the fact that $c_\gamma$ in (\ref{eik567}) changes 
its sign at $t=1$. It is even possible to build an autonomous counter-example up 
to start with a front with several connected components;\\ 
(ii) $c_0 \geq 0$, and 
therefore it complements also the results of Section~\ref{slep};  
indeed the unique solution $u$ of (\ref{slepcev-intro}) has the same 0--level-set than $U_1$ 
(obtained with $\gamma(t)\equiv 1$) and, with the notations of Theorem \ref{th:4},  
$\rho^+= \1_{\{U^1(\cdot ,t)\geq 0\}}$ and  
$\rho^-= \1_{\{U^1(\cdot ,t)> 0\}}.$ In particular, for $t\geq 1,$ 
$$ 
\{ v^+(\cdot,t)\geq 0\}={\{U^1(\cdot ,t)\geq 0\}}=[-y_1(t),y_1(t)] 
$$ 
and 
$$ 
 \{ v^-(\cdot,t)> 0\}={\{U^1(\cdot ,t)> 0\}}=\emptyset. 
$$ 
Finally, we note that there is no strong solutions since (\ref{FormeFaible3}) is obviously never satisfied; 
\\ 
(iii) $c_0=1$ 
does not satisfies {\bf (H1)} but because of the finite speed of propagation property, 
it is possible to keep the same 
solution on a large ball in space and for $t\in (0,T)$, if we replace $c_0$ by 
a function with compact support in space such that $c_0(x,t)=1$ for 
$|x|\le R$ with $R$ large enough. By this way, it is possible for $c_0$ to 
satisfy {\bf (H1)}.  
\end{Example} 
 
\section{Uniqueness results for weak solutions of  (\ref{dislocation})}\label{sect4} 
 
Uniqueness of weak solutions of  (\ref{dislocation}) is false in general as 
shown in the counter-example of the previous section for sign changing 
velocities $c_1$. This is in particular related to the  
``fattening phenomenon".  
In \cite{acm05} and \cite{bl06}  
the authors proved that there is a unique ``classical" viscosity solution for (\ref{dislocation}) under the assumptions that the initial set $\{u(\cdot,0)\geq 0\}$ 
has the ``interior sphere property" and that $c_1(\cdot,t) \geq |c_0(\cdot,t)|_{L^1}$ for any $t\geq 0$---condition which ensures that the velocity $\bar{c}$ is nonnegative.  
By ``classical" continuous viscosity solutions we mean that $t\to \1_{\{ u(\cdot,t)\geq 0\}}$ is continuous in $L^1$, which entails that $(x,t)\mapsto \bar{c}(x,t)$ is continuous, and that 
(\ref{dislocation}) holds in the usual viscosity sense.  
 
Here we prove Theorem \ref{th:2}. If the condition  $c_1(\cdot,t) \geq |c_0(\cdot,t)|_{L^1}$ is satisfied, then weak solutions are viscosity solutions. 
We also prove that the weak solution is unique if we suppose moreover either that the initial condition has the interior sphere property, or 
that the strict inequality $c_1(\cdot,t) > |c_0(\cdot,t)|_{L^1}$ holds for any 
$t\geq 0$. \\ 
 
\noindent{\bf Proof of Theorem \ref{th:2}.} \\ 
\noindent{\it 1. Weak solutions are classical continuous viscosity solutions.} 
Let $u$ be a weak solution and let ${\bar c}$ be associated with $u$ as in 
Definition \ref{defi:1}.  
Then, for any $x\in \R^N$ and for almost all $t\in [0,T]$, we have 
$$ 
\begin{array}{rl} 
{\bar c}(x,t) \;  \geq & \displaystyle{  c_1(x,t) + c^+_0 (\cdot,t) \star \1_{\{u(\cdot,t) > 0\}}(x) - c^-_0(\cdot,t) \star \1_{\{u(\cdot,t)\geq  
0\}}(x)  }  \\ 
 & \\ 
\geq & \displaystyle{ c_1(x,t)- |c_0^- (\cdot, t)|_{L^1}  }\\ 
 & \\ 
\geq &  \delta \;\geq\;0\;. 
\end{array} 
$$ 
From \cite[Theorem 4.2]{ley01}, there exists a constant $\eta$ 
which depends on $T$ such that (\ref{lgb1111}) implies 
\begin{equation}\label{lbgT} 
- |u| - |D u| \leq -\eta \quad\hbox{on  } \R^N\times (0,T) \;, 
\end{equation} 
when we assume moreover that ${\bar{c}}$ is continuous. In our case where 
${\bar{c}}$ is not assumed continuous in time, (\ref{lgb1111}) follows from the 
$L^1$-stability result Theorem \ref{L1stability} where we approximate 
${\bar{c}}$ by a continuous function, and from the usual stability for 
$L^1$-viscosity subsolutions.\\ 
Let us note that from the proof of \cite[Corollary 2.5]{bl06}, we 
have in the viscosity sense 
\begin{equation}\label{lbgTx} 
- |u(\cdot,t)| - |D u(\cdot,t)| \leq -\eta \quad\hbox{on  } \R^N,\ 
\mbox{for every}\ t\in (0,T) \; . 
\end{equation} 
Following  
\cite[Corollary 2.5]{bl06}, we get that, for every $t\in (0,T)$, the 
0--level-set of $u(\cdot,t)$ has a  zero Lebesgue measure.  
Then we deduce that  
$$ 
\chi(x,t)=\1_{\left\{u(\cdot,t)\ge 0\right\}}(x), \quad\mbox{for 
  a.e.}\ x\in\R^N,\ \mbox{for all} \ t\in (0,T) 
$$ 
which (with (\ref{FormeFaible2})) 
entails that  
$$ 
{\bar c}(x,t)=c_1 + c_0 \star \1_{\{u(\cdot,t) \geq 0\}} 
$$ 
for any $(x,t)$. Moreover $t\mapsto \1_{\{u(\cdot,t) \geq 0\}}$ is also 
continuous in $L^1$, and then ${\bar{c}}$ is continuous. Therefore $u$ is a 
classical viscosity solution of (\ref{dislocation}).\\ 
 
\noindent{\it 2. Uniqueness when $u_0$ is semiconvex (Part (ii)).} If we assume that (\ref{HypPos}) and (\ref{lgb1111}) hold and that $u_0$ is semiconvex,  
then weak solutions are viscosity solutions, and we can apply the 
uniqueness result for viscosity solutions given in \cite{bl06}, namely 
Theorem 4.2 (which remains true under our assumptions) which requires in particular semiconvexity of the velocity, 
see assumption {\bf (H2)}. \\ 
 
\noindent{\it 3. A Gronwall type inequality (Part (i)).} From now on we assume that $\delta>0$ and we aim at proving that the  
solution to (\ref{dislocation}) is unique. Let $u_1$, $u_2$ be two solutions. We set  
$$ 
\rho_i= \1_{\{u_i(\cdot,t) \geq 0\}}\qquad {\rm and }\qquad  
{\bar{c}}_i(x,t)=c_0\star  \rho_i+c_1\qquad \mbox{\rm for $i=1,2$.} 
$$ 
We want to prove in a first step the following Gronwall type inequality for any $t$ sufficiently small: 
\begin{equation}\label{EstiPeri} 
\begin{array}{l} 
|\rho_1(\cdot,t)-\rho_2(\cdot,t)|_{L^1}\\ 
\; \leq \; \displaystyle{   C \left[  {\rm per}(\{u_1(\cdot,t)\geq 0\})+{\rm per}(\{u_2(\cdot,t)\geq 0\})\right] \int_0^t|\rho_1(\cdot,s)-\rho_2(\cdot,s)|_{L^1}ds   } 
\end{array} 
\end{equation} 
where $C$ is a constant depending on the constants of the problem, where ${\rm per}(\{u_i(\cdot,t)\geq 0\})$ is the $\mathcal{H}^{N-1}$ measure of  
the set $\partial \{u_i(\cdot,t)\geq 0\})$ (for $i=1,2$).  
 
We have 
\begin{equation}\label{toti} 
|\rho_1(\cdot,t)-\rho_2(\cdot,t)|_{L^1}\leq \mathcal{L}^N (\{-\alpha_t\leq u_1(\cdot,t)< 0\})+\mathcal{L}^N (\{-\alpha_t\leq u_2(\cdot,t)< 0\})\;, 
\end{equation} 
where $\mathcal{L}^N$ is the Lebesgue measure in $\R^N$ and 
\begin{equation}\label{defalpha} 
\alpha_t =\sup_{s\in [0,t]}|(u_1-u_2)(\cdot,s)|_\infty\qquad \mbox{\rm for any $t\in (0,T)$.} 
\end{equation} 
In order to estimate the right-hand side of inequality (\ref{toti}),  
as in the proof of \cite[Theorem 4.2]{bl06},  
we need a lower-gradient bound as well 
as a semiconvexity property for $u_1$ and $u_2$. We already know from step 
1 that ${\bar{c}}_i$ is continuous for $i=1,2$. 
  
Let us start to estimate the right-hand side of (\ref{toti}). From the 
``stability estimates'' on the solutions with respect to variations of the velocity (see  \cite[Lemma 2.2]{bl06}), we have 
\begin{equation}\label{eq:stab0} 
\alpha_t\leq |Du_0|_\infty e^{L t}\int_0^t |({\bar{c}}_1-{\bar{c}}_2)(\cdot,s)|_{\infty }ds\;,    
\end{equation} 
where $L=L_0+L_1$. 
Therefore 
\begin{equation}\label{EstiDeltat} 
\alpha_t\leq m_0 \; |Du_0|_\infty e^{L t}\int_0^t |(\rho_1-\rho_2)(\cdot,s)|_{L^1}ds\;. 
\end{equation} 
where the constant $m_0$ is given in {\bf (H2)}.  
In particular, since the $\rho_i(\cdot,t)$ are continuous in $L^1$ and equal at time $t=0$ for $i=1,2$, we have $\alpha_t/t\to 0$ as $t\to 0^+$.  
 
From now on, we mimick the proof of \cite[Proposition 4.5]{bl06}.  
Using 
the lower-gradient bound (\ref{lbgT}) for $u_i(\cdot,t)$ combined with 
the increase principle (see \cite[Lemma 2.3]{bl06}), we obtain for $\alpha_t < \eta/2$ that  
$$ 
\{-\alpha_t\leq u_i(\cdot,t)< 0\}\subset \{ u_i(\cdot,t)\geq 0\}+(2\alpha_t/\eta) \overline{B}(0,1)\;, 
$$ 
for $i=1,2$.  
From the interior ball regularization Lemma \ref{lem:bornebouleint}, 
the set $\{ u_i(\cdot,t)\geq 0\}$ satisfies for $t\in (0,T)$ the interior 
ball property of radius $r_t=\eta t/C_0.$ 
Applying 
\cite[Lemma 2.5 and 2.6]{acm05}, we obtain for $\sigma_t=2\alpha_t/\eta$ that  
$$ 
\begin{array}{rl} 
\mathcal{L}^N(\{-\alpha_t \leq u_i(\cdot,t)< 0\})  \leq \; & \displaystyle{         
\mathcal{L}^N\left((\{ u_i(\cdot,t)\geq 0\}+\sigma_t \, \overline{B}(0,1))\backslash \{ u_i(\cdot,t)\geq 0\}\right)        }\\ 
 & \\ 
\leq  &  \displaystyle{   \frac{r_t }{N}\left[ \left(1+ \frac{\sigma_t}{r_t}\right)^N-1\right]  {\rm per }(\{ u_i(\cdot,t)\geq 0\})    }\\ 
& \\ 
\leq &  \displaystyle{    \frac{  2^{N}\alpha_t}{\eta}  {\rm per }(\{ u_i(\cdot,t)\geq 0\})    } 
\end{array}$$ 
(using $(1+a)^N-1 \le a N (1+a)^{N-1}$ for $a\ge 0$) for $t\in [0,\tau]$ where 
$0<\tau\leq T$ is defined by  
\begin{eqnarray} \label{deftau} 
\tau = \mathop{\rm sup} 
\{ t>0 : \alpha_t <\frac{\eta}{2} \ {\rm and} \ 2\frac{C_0}{\eta^2} 
 \frac{\alpha_t}{t} \le 1 \}. 
\end{eqnarray} 
 Putting together (\ref{EstiDeltat}), (\ref{toti}) and the previous inequality proves (\ref{EstiPeri}).\\ 
 
\noindent{\it 4. Uniqueness when $\delta>0$ (Part (i)).} We now complete the uniqueness proof under the assumption $\delta>0$. 
For this we first show that $\rho_1=\rho_2$ in $[0,\tau].$ In order to apply Gronwall Lemma to the $L^1$-estimate (\ref{EstiPeri})  
obtained in Step 3, it is enough to prove that the functions $t\mapsto  {\rm per}(\{u_i(\cdot,t)\geq 0\})$ 
belong to $L^1$. For this, let us set 
$$ 
w_i(x)=\inf\{t\geq 0\; :\; u_i(x,t)\geq 0\}\;. 
$$ 
Since $u_i$ solves the eikonal equation $(u_i)_t= 
{\bar{c}}_i(x,t)|Du_i|$, from classical representation formulae, we have
\begin{eqnarray*} 
 \{u_i(\cdot ,t)\geq 0\} &=& \{x \, : \ \exists y(\cdot),\ |\dot{y}(s)|\leq
c(y(s),s), \  0\leq s\leq t, \\ 
&& \hspace*{3.5cm} u_0(y(0))\geq 0 \ {\rm and} \ y(t)=x \}.
\end{eqnarray*} 
Therefore
\begin{eqnarray*} 
w_i(x)&=& {\rm inf} \{ t\geq 0\; :\;  \ \exists y(\cdot),\ |\dot{y}(s)|\leq
c(y(s),s), \   0\leq s\leq t, \\
&& \hspace*{4.2cm} u_0(y(0))\geq 0 \ {\rm and} \ y(t)=x \}.
\end{eqnarray*}
Applying the dynamic programming principle, since ${\bar{c}}_i\geq  \delta>0,$
we obtain that $w_i$ is Lipschitz
continuous and is a viscosity solution of the autonomous equation 
${\bar{c}}_i(x,w_i(x))|Dw_i(x)|=1$.  
Note  that $\{u_i(\cdot,t)\geq 0\}=\{w_i\leq t\}$. In particular, by Theorem \ref{basicthm} (i), 
$\{w_i\leq t\}\subset B(0,R_0+(M_0+M_1)T)$ is bounded for any $t.$ From the 
coarea formula, we have 
$$ 
\begin{array}{rl} 
\displaystyle{  \int_0^t {\rm per }(\{u_i(\cdot,s)\geq 0\})ds\; } = & \displaystyle{   \int_0^t {\rm per }(\{w_i\leq s\})ds  }\\ 
 & \\ 
 = & \displaystyle{  \int_{\{w_i\leq t\}} |Dw_i(x)|dx \;,} 
 \end{array} 
$$ 
which is finite since $w_i$ is Lipschitz continuous. Therefore we have proved that $t\mapsto  {\rm per}(\{u_i(\cdot,t)\geq 0\})$ 
belongs to $L^1([0,\tau])$, which entails from Gronwall Lemma that 
$\rho_1=\rho_2$ in $[0,\tau]$ 
since $\rho_1(\cdot,0)=\rho_2(\cdot,0)$.  
Hence ${\bar{c}}_1={\bar{c}}_2$ and $u_1=u_2$ in  $[0,\tau].$ From the 
definition of $\alpha_t$ and $\tau,$ (see (\ref{defalpha}) and (\ref{deftau})), necessarily  $\tau=T.$ 
It completes the proof. 
\hfill $\Box$

\section{Nonnegative kernel $c_0$ and Slep$\rm \check{c}$ev formulation for 
the nonlocal term}\label{slep}

In this section, we deal with nonnegative kernels $c_0\ge 0$. In this monotone 
framework, inclusion principle for evolving sets and comparison for solutions 
to the dislocation equation are expected (see Cardaliaguet \cite{cardaliaguet00}  
for related results). 
We start by studying the right level-set equation using a Slep$\rm 
\check{c}$ev formulation with the convolution  term using all the level-sets 
$\left\{u(\cdot,t)\ge u(x,t)\right\}$ instead of only one level-set 
$\left\{u(\cdot,t)\ge 0\right\}$. This choice is motivated by the good 
stability properties of the Slep$\rm \check{c}$ev formulation.\\ 
The equation we are concerned with is 
\begin{equation}\label{slepcev} 
\left\{\begin{array}{ll} 
\displaystyle{\frac{\partial u}{\partial t}   
= c^{+}[u](x,t)|D u| } 
& \mbox{in} \ \R^N\times (0,T)\\ 
u(\cdot,0)=u_0 &  \mbox{in} \ \R^N\; , 
\end{array}\right. 
\end{equation} 
where the nonlocal velocity is 
\begin{equation}\label{eq:c>=} 
\begin{array}{lcl} 
c^{+}[u](x,t) &=&  c_1(x,t) + c_0(\cdot,t) \star \1_{\{u(\cdot,t)\geq 
  u(x,t)\}}(x)  \\ 
\\ 
&=& \displaystyle{c_1(x,t) + \int_{\R^N} c_0(x-z,t)\1_{\{u(\cdot,t)\geq u(x,t)\}} (z)dz } 
\end{array}   
\end{equation} 
and the additional velocity $c_1$ has no particular sign.  
 
We denote 
\begin{eqnarray*} 
c^{-}[u](x,t)= c_1(x,t) + \int_{\R^N} c_0(x-z,t)\1_{\{u(\cdot,t)> u(x,t)\}} (z)dz . 
\end{eqnarray*} 
 
We recall the notion of viscosity solutions for (\ref{slepcev}) as it 
appears in  \cite{slepcev03}.  
 
\begin{Definition}\label{Defi-slepcev}{\bf (Slep$\rm \check{c}$ev viscosity 
 solutions)}\\ 
An upper-semicontinuous function $u:\R^N\times [0,T]\to \R$ is a viscosity  
subsolution of 
(\ref{slepcev}) if, for any $\varphi\in C^1(\R^N\times [0,T]),$ for any maximum  
point 
$(\bar{x},\bar{t})$ of $u-\varphi,$ if $\bar{t}>0$ then 
$$ 
\frac{\partial \varphi}{\partial t}(\bar{x},\bar{t})\leq  
c^{+}[u](\bar{x},\bar{t})|D\varphi (\bar{x},\bar{t})|  
$$ 
and $u(\bar{x},0)\leq u_0 (\bar{x})$ if $\bar{t}=0.$ 
 
A lower-semicontinuous function $u:\R^N\times [0,T]\to \R$ is a viscosity  
supersolution of 
(\ref{slepcev}) if, for any $\varphi\in C^1(\R^N\times [0,T]),$ for any minimum  
point 
$(\bar{x},\bar{t})$ of $u-\varphi,$ if $\bar{t}>0$ then 
$$ 
\frac{\partial \varphi}{\partial t}(\bar{x},\bar{t})\geq  
c^{-}[u](\bar{x},\bar{t})|D\varphi (\bar{x},\bar{t})|  
$$ 
and $u(\bar{x},0)\geq u_0 (\bar{x})$ if $\bar{t}=0.$ 
 
A locally bounded function is a viscosity solution of (\ref{slepcev}) if its 
upper-semicontinuous enveloppe is subsolution and its lower-semicontinuous 
envelope is supersolution of (\ref{slepcev}). 
\end{Definition} 
 
Note that for the 
supersolution, we require the viscosity inequality with $c^{-}$ instead of $c^{+}$. It is the definition providing the expected stability results  
(see \cite{slepcev03}). 
 
\begin{Theorem} \label{comp-slep}{\bf (Comparison principle)}\\ 
Assume {\bf (H0')}, and that the kernel $c_0\ge 0$ and $c_1$ satisfy  
{\bf (H1)}.  
Let $u$ (respectively $v$)  
be a bounded upper-semicontinuous subsolution (respectively a bounded lower semicontinuous supersolution) 
of (\ref{slepcev}). 
Then $u\leq v$ in $\R^N\times [0,T].$ 
\end{Theorem} 
 
\begin{Remark} \label{rem-slep1} \rm We could deal with second-order terms in (\ref{slepcev}) (for instance 
we can add the mean curvature to the velocity (\ref{eq:1})). 
See Forcadel \cite{forcadel05}  and Srour \cite{srour06} for related results. 
\end{Remark} 
 
Before giving the proof of Theorem \ref{comp-slep}, let us note the 
following consequence. \\ 
 
\noindent {\bf Proof of Theorem \ref{th:3}.} 
The uniqueness of a continuous viscosity solution to (\ref{slepcev}) is  
an immediate consequence of Theorem \ref{comp-slep}. Then existence is proved 
by Perron's method using classical arguments (see for instance \cite[Theorem 1.2]{dlks04}), so we 
skip the details. 
\hfill$\Box$ \\ 
 
\noindent{\bf Proof of Theorem \ref{comp-slep}.} \\ 
\noindent{\it 1. The test-function.} 
 Since $u-v$ is a bounded upper-semicontinuous function,   
for any $\varepsilon, \eta, \alpha >0$ and $K=2(L_0+L_1)\ge 0,$ the supremum 
\begin{eqnarray*}  
M_{\varepsilon, \eta, \alpha}= \mathop{\rm sup}_{(x,y,t)\in (\R^N)^2\times [0,T]} 
\{ u(x,t)-v(y,t)-{\rm e}^{Kt}(\frac{|x-y|^2}{\varepsilon^2}+\alpha |x|^2 +\alpha |y|^2)-\eta t\} 
\end{eqnarray*} 
is finite and achieved at a point $(\bar{x},\bar{y},\bar{t}).$ 
 Classical arguments show that 
 \begin{eqnarray*} 
 \liminf_{\varepsilon, \eta, 
 \alpha \to 0}M_{\varepsilon, \eta, \alpha}  =   
 \mathop{\rm sup}_{\R^N\times [0,T]} \{ u-v\} 
 \end{eqnarray*} 
and that 
\begin{eqnarray} \label{borne123} 
\frac{|\bar{x}-\bar{y}|^2}{\varepsilon^2}, \alpha |\bar{x}|^2, \alpha |\bar{y}|^2 \leq M_\infty, 
\end{eqnarray} 
where $M_\infty=|u|_\infty +|v|_\infty.$ \\ 
 
\noindent{\it 2. Viscosity inequalities when $\bar{t}>0.$} 
Writing the viscosity inequalities for the subsolution $u$ and the supersolution $v,$ 
we obtain  
\begin{equation}\label{ineq333} 
\hspace*{0.9cm} K{\rm e}^{K\bar{t}}(\frac{|\bar{x}-\bar{y}|^2}{\varepsilon^2}+ 
 \alpha |\bar{x}|^2+ \alpha |\bar{y}|^2  ) +\eta   
\leq c^{+}[u](\bar{x},\bar{t})|\bar{p}+\bar{q}_x|-c^{-}[v](\bar{y},\bar{t})|\bar{p}-\bar{q}_y|  
\end{equation} 
where $\bar{p} = 2e^{K\bar{t}}{(\bar{x}-\bar{y})}/{\varepsilon^2}$, $\bar{q}_x=2e^{K\bar{t}}\alpha \bar{x}$ and 
$\bar{q}_y=2e^{K\bar{t}}\alpha \bar{y}.$ 
We point out a difficulty to obtain this inequality: in general, one gets 
it by doubling the time variable first and then by passing to the limit in the 
time penalization. This is not straightforward here because of the dependence 
with respect to time of the nonlocal terms. But the stability arguments of 
the Slep$\rm\check{c}$ev formulation take care of this difficulty. 
\\ 
 
\noindent{\it 3. Difference between $\{u(\cdot ,\bar{t})\geq u(\bar{x} ,\bar{t})\}$ 
and $\{v(\cdot ,\bar{t})> v(\bar{y} ,\bar{t})\}$.} We have 
\begin{eqnarray} \label{inclu1} 
\{u(\cdot ,\bar{t})\geq u(\bar{x} ,\bar{t})\}\subset \{v(\cdot ,\bar{t})> v(\bar{y} ,\bar{t})\}\cup 
\mathcal{E}, 
\end{eqnarray} 
where $\mathcal{E}=  \{u(\cdot ,\bar{t}) \geq u(\bar{x} ,\bar{t})\} \cap 
\{v(\cdot ,\bar{t})\leq  v(\bar{y} ,\bar{t})\}.$ If $x\in  \mathcal{E},$ 
then $u(\bar{x} ,\bar{t})-v(\bar{y} ,\bar{t})\leq u(x ,\bar{t})-v(x ,\bar{t}).$ 
But from the definition of $M_{\varepsilon, \eta, \alpha},$ 
\begin{eqnarray*} 
&& u(x ,\bar{t})-v(x ,\bar{t})-{\rm e}^{K\bar{t}}2\alpha |x|^2-\eta \bar{t} \\ 
&\leq& u(\bar{x} ,\bar{t})-v(\bar{y} ,\bar{t}) 
-{\rm e}^{K\bar{t}}(\frac{|\bar{x}-\bar{y}|^2}{\varepsilon^2} + \alpha |\bar{x}|^2+ \alpha |\bar{y}|^2) -\eta \bar{t}. 
\end{eqnarray*} 
It follows that  
\begin{eqnarray*} 
 \mathcal{E} \subset \{x\in \R^N : 
|x|^2\geq \frac{1}{2} ( |\bar{x}|^2+ |\bar{y}|^2) + \frac{|\bar{x}-\bar{y}|^2}{2\alpha \varepsilon^2}\}. 
\end{eqnarray*} 
\smallskip 
 
\noindent{\it 4. Upper bound for $c^{+}[u](\bar{x},\bar{t}).$} 
We have 
\begin{eqnarray} \label{ineq222} 
\hspace*{0.3cm} c^{+}[u](\bar{x},\bar{t}) &=& 
\int_{\R^N} c_0(\bar{x}-z,\bar{t}) \1_{ \{u(\cdot ,\bar{t})\geq u(\bar{x} ,\bar{t})\}}(z)dz +c_1(\bar{x},\bar{t}) \\ 
&\leq & 
\int_{\R^N} (c_0(\bar{x}-z,\bar{t})-c_0(\bar{y}-z,\bar{t})) \1_{ \{u(\cdot ,\bar{t})\geq u(\bar{x} ,\bar{t})\}}(z)dz  
\nonumber\\ 
&& + \int_{\R^N} c_0(\bar{y}-z,\bar{t}) \1_{ \{u(\cdot ,\bar{t})\geq u(\bar{x} ,\bar{t})\}}(z)dz 
+c_1(\bar{x},\bar{t}). 
\nonumber 
\end{eqnarray} 
Using that $c_0\geq 0$ and (\ref{inclu1}), we obtain 
\begin{eqnarray*} 
 \int_{\R^N} c_0(\bar{y}-z,\bar{t}) \1_{ \{u(\cdot ,\bar{t})\geq u(\bar{x} ,\bar{t})\}}(z)dz 
\leq \int_{ \{v(\cdot ,\bar{t})> v(\bar{y} ,\bar{t})\}\cup 
\mathcal{E}  } c_0(\bar{y}-z,\bar{t})dz. 
\end{eqnarray*} 
From (\ref{ineq222}), we get 
\begin{equation}\label{eq:I1I2} 
c^{+}[u](\bar{x},\bar{t}) \quad \leq \quad  
c^{-}[v](\bar{y},\bar{t}) + \mathcal{I}_1 + \mathcal{I}_2 +c_1(\bar{x},\bar{t})-c_1(\bar{y},\bar{t}),    
\end{equation} 
where 
\begin{eqnarray*} 
\mathcal{I}_1= 
\int_{\R^N} (c_0(\bar{x}-z,\bar{t})-c_0(\bar{y}-z,\bar{t})) \1_{ \{u(\cdot ,\bar{t})\geq u(\bar{x} ,\bar{t})\}}(z)dz  
\end{eqnarray*} 
and 
\begin{eqnarray*} 
\mathcal{I}_2=  \int_{\mathcal{E}} c_0(\bar{y}-z,\bar{t})dz. 
\end{eqnarray*} 
 
\noindent{\it 5. Estimate of $\mathcal{I}_1$ using {\bf (H1)}.} 
We have 
\begin{eqnarray*} 
c_0(\bar{x}-z,t)-c_0(\bar{y}-z,t)=\int_0^1 D_x c_0((1-\lambda)(\bar{y}-z)+\lambda(\bar{x}-z), \bar{t}) 
(\bar{x}-\bar{y})d\lambda. 
\end{eqnarray*} 
It follows 
\begin{eqnarray} 
|\mathcal{I}_1|&\leq& \int_{\R^N}\int_0^1 |D_x c_0((1-\lambda)(\bar{y}-z)+\lambda(\bar{x}-z), \bar{t})| 
|\bar{x}-\bar{y}|d\lambda dz \nonumber \\ 
&\leq& |D_x c_0(\cdot , \bar{t})|_{L^1}|\bar{x}-\bar{y}| \nonumber\\ 
&\leq & L_0|\bar{x}-\bar{y}|. \label{1ereinteg} 
\end{eqnarray} 
 
\noindent{\it 6. Estimate of the right-hand side of inequality (\ref{ineq333}).} 
Noticing that   
$|\bar{p}||\bar{x}-\bar{y}|=2 {\rm e}^{K\bar{t}} |\bar{x}-\bar{y}|^2/\varepsilon^2$ and using 
(\ref{1ereinteg}), (\ref{eq:I1I2}) and {\bf (H1)}, 
we have 
\begin{eqnarray*} 
& & c^{+}[u](\bar{x},\bar{t})|\bar{p}+\bar{q}_x|-c^{-}[v](\bar{y},\bar{t})|\bar{p}-\bar{q}_y| \\ 
& \leq &  
(c^{-}[v](\bar{y},\bar{t}) + \mathcal{I}_1 + \mathcal{I}_2 +c_1(\bar{x},\bar{t})-c_1(\bar{y},\bar{t})) 
|\bar{p}+\bar{q}_x|  
 - c^{-}[v](\bar{y},\bar{t})|\bar{p}-\bar{q}_y| \\ 
&\leq& 
 \left|c^{-}[v](\bar{y},\bar{t})\right| |\bar{q}_x +\bar{q}_y | 
+(L_0+L_1)|\bar{x}-\bar{y}||\bar{p}+\bar{q}_x| +\mathcal{I}_2|\bar{p}+\bar{q}_x| \\ 
&\leq&  
 (M_0+M_1) (|\bar{q}_x | + | \bar{q}_y |) 
+ 2 {\rm e}^{K\bar{t}} (L_0+L_1)\frac{|\bar{x}-\bar{y}|^2}{\varepsilon^2} \\ 
&& 
+(L_0+L_1)|\bar{x}-\bar{y}| |\bar{q}_x|  
+ \mathcal{I}_2 (|\bar{q}_x|+ 2{\rm e}^{K\bar{t}}\frac{|\bar{x}-\bar{y}|}{\varepsilon^2}). 
\end{eqnarray*} 
Since $|\bar{q}_x|, |\bar{q}_y| \to 0$ as $\alpha\to 0$ (see (\ref{borne123})) and 
$\mathcal{I}_2$ is bounded by $|c_0(\cdot ,\bar{t})|_{L^1}\leq L_0$, 
there exists a modulus $m_\varepsilon (\alpha)\to 0$ as $\alpha\to 0$ such that 
(\ref{ineq333}) becomes 
\begin{eqnarray*} 
&& K{\rm e}^{K\bar{t}}(\frac{|\bar{x}-\bar{y}|^2}{\varepsilon^2}+ 
 \alpha |\bar{x}|^2+ \alpha |\bar{y}|^2  ) +\eta \\ 
&\leq&   
m_\varepsilon (\alpha) 
+2(L_0+L_1)  {\rm e}^{K\bar{t}}   \frac{|\bar{x}-\bar{y}|^2}{\varepsilon^2} 
+2 \, \mathcal{I}_2 \, {\rm e}^{K\bar{t}}\frac{|\bar{x}-\bar{y}|}{\varepsilon^2}. 
\end{eqnarray*} 
Recalling that we chose $K\geq 2(L_0+L_1),$ 
we finally obtain 
\begin{eqnarray} \label{avantlimite} 
0< \eta\leq m_\varepsilon (\alpha)+2 \, \mathcal{I}_2 \, {\rm e}^{K\bar{t}}\frac{|\bar{x}-\bar{y}|}{\varepsilon^2}. 
\end{eqnarray} 
 
\noindent{\it 7. Limit when $\alpha\to 0.$} 
First, suppose that 
\begin{eqnarray} \label{premiercas1} 
\frac{|\bar{x}-\bar{y}|^2}{\varepsilon^2} \to 0 \ \ \ {\rm as} \ \alpha\to 0. 
\end{eqnarray} 
It follows that $|\bar{x}-\bar{y}|\to 0$ as $\alpha\to 0.$ Passing to the limit in  
(\ref{avantlimite}), we obtain a contradiction. 
Therefore, (\ref{premiercas1}) cannot hold and, up to extract a subsequence, 
there exists $\delta >0$ such that 
\begin{eqnarray} \label{casdelta} 
\frac{|\bar{x}-\bar{y}|^2}{\varepsilon^2} \geq \delta > 0 \ \ \ {\rm for} \ \alpha >0 \ {\rm small \ enough}. 
\end{eqnarray} 
From (\ref{avantlimite}) and (\ref{borne123}), we get 
\begin{eqnarray} \label{abcde} 
\eta \leq \mathop{\rm lim\,sup}_{\alpha\to 0} 2 \, \mathcal{I}_2 \,  
{\rm e}^{K\bar{t}}\frac{|\bar{x}-\bar{y}|}{\varepsilon^2} 
\leq \frac{2 {\rm e}^{K\bar{t}} M_\infty^{1/2}}{{\varepsilon}} \mathop{\rm lim\,sup}_{\alpha\to 0}\mathcal{I}_2. 
\end{eqnarray} 
To obtain a contradiction, it suffices to show that ${\rm lim\,sup}_{\alpha\to 0}\mathcal{I}_2 = 0.$  
 
\noindent{\it 8. Convergence of $\mathcal{I}_2$ to $0$ when $\alpha\to 0.$} 
By a change of variable, we have 
\begin{eqnarray*} 
\mathcal{I}_2=  \int_{\mathcal{E}} c_0(\bar{y}-z,\bar{t})dz 
\leq  \int_{\bar{\mathcal{E}}} c_0(z,\bar{t})dz, 
\end{eqnarray*} 
where 
\begin{eqnarray*} 
\bar{ \mathcal{E}} = \{x\in \R^N : 
|x-\bar{y}|^2\geq \frac{1}{2} ( |\bar{x}|^2+ |\bar{y}|^2)  
+ \frac{|\bar{x}-\bar{y}|^2}{2\alpha \varepsilon^2}\}. 
\end{eqnarray*} 
Since $|c_0(\cdot ,\bar{t})|_{L^1}\leq L_0$, 
to prove that $\mathcal{I}_2 \to 0,$ it suffices to show that $\bar{ \mathcal{E}}\subset \R^N\backslash  
B(0,R_\alpha)$ with $R_\alpha\to +\infty.$ From (\ref{casdelta}), if $x\in 
\bar{ \mathcal{E}},$ then 
\begin{eqnarray*} 
|x|^2&\geq& -2|x| |\bar{y}|+\frac{1}{2} |\bar{x}|^2 -\frac{1}{2} |\bar{y}|^2 +\frac{\delta}{2\alpha} \\ 
&\geq & 
-2|x| |\bar{y}| - |\bar{y}||\bar{x}-\bar{y}|+\frac{\delta}{2\alpha} \\ 
&\geq & 
\frac{1}{2\alpha} (\delta - 2(C^2+2C|x|)\sqrt{\alpha}), 
\end{eqnarray*} 
since by (\ref{borne123}), there exists $C>0$ such that $|\bar{x}|, |\bar{y}|\leq C/\sqrt{\alpha}$ 
and $|\bar{x}-\bar{y}|\leq C.$ It follows that 
\begin{eqnarray*} 
|x|\geq \frac{1}{\sqrt{\alpha}}\left( -C+\sqrt{C^2+\delta/2 -C^2\sqrt{\alpha}}\right) := R_\alpha 
\mathop{\longrightarrow}_{\alpha\to 0} +\infty. 
\end{eqnarray*} 
 
\noindent{\it 9. End of the proof.} 
Finally, for every $\varepsilon,$ 
if $\alpha=\alpha_\varepsilon$ is small enough, the supremum $M_{\varepsilon, \eta, \alpha}$  
is necessarily achieved for $\bar{t}=0.$ It follows 
$$ 
M_{\varepsilon, \eta, \alpha}\leq u(\bar{x},0)-v(\bar{y},0)- \frac{|\bar{x}-\bar{y}|^2}{\varepsilon^2} 
\leq u_0(\bar{x},0)-u_0(\bar{y},0)- \frac{|\bar{x}-\bar{y}|^2}{\varepsilon^2}. 
$$ 
Since, $u_0$ is uniformly continuous, for all $\rho >0,$ there exists $C_\rho>0$ such that 
$$ 
M_{\varepsilon, \eta, \alpha}\leq \rho + C_\rho |\bar{x}-\bar{y}|-\frac{|\bar{x}-\bar{y}|^2}{\varepsilon^2} 
\leq \rho + \frac{C_\rho^2 \varepsilon^2}{4}. 
$$ 
Passing to the limits $\varepsilon\to 0$ and then $\rho, \alpha, \eta \to 0,$ we 
obtain that ${\rm sup}\{u-v\}\leq 0$ 
\hfill $\Box$ \\


Now we turn to the connections with \emph{discontinuous solutions} and weak solutions, which are closely connected.  
To do so, if $u$ is the unique continuous solution of (\ref{slepcev}) given by 
Theorem~\ref{th:3}, we recall that we use the notations  
\begin{eqnarray*} 
\rho^+ : =\1_{\{u\geq 0\}}, \  \ \rho^- : =\1_{\{u > 0\}} \  \ {\rm and} \ \   
c[\rho](x,t) = c_0 (\cdot ,t) \star \rho(\cdot,t) (x) +c_1(x,t). 
\end{eqnarray*} 
 
\noindent{\bf Proof of Theorem~\ref{th:4}.} \\ 
{\it 1. Claim: Under the assumptions of Theorem~\ref{th:4}, the functions 
 $\rho^+$ and $\rho^-$ are  $L^1$-viscosity solutions of the equation} 
\begin{eqnarray} \label{eik111} 
\left\{ 
\begin{array}{ll} 
\rho_t = c[\rho]|D\rho| & {in  } \ \R^N\times [0,T)\; , \\[2mm] 
\rho(x,0) = \1_{\{u_0  \geq 0\}} & {in  } \  \R^N\; . 
\end{array} 
\right. 
\end{eqnarray} 
We consider two sequences of smooth nondecreasing functions $(\psi_\alpha)_\alpha$, $(\psi^\alpha)_\alpha$, taking  
values in $[0,1]$, such that, for any $s\in \R$ 
$$  
\psi_\alpha (s) \leq \1_{\{ x > 0\}} (s) \leq \1_{\{x \geq 0\}} (s) \leq \psi^\alpha (s) \; , 
$$ 
and such that, as $\alpha \to 0$, $ \psi_\alpha \uparrow \1_{\{x > 0\}}$, $\psi^\alpha \downarrow \1_{\{x \geq 0\}}$. 
 
We first remark that $u$ satisfies, in the sense of Definition \ref{Defi-slepcev}, 
\begin{eqnarray*} 
u_t &\leq & c[\psi^\alpha (u(\cdot,t)-u(x,t))]|Du| \quad \mbox{in} \ \R^N\times (0,T)\; , \\ 
u_t & \geq & c[\psi_\alpha (u(\cdot,t)-u(x,t))]|Du| \quad \mbox{in} \  \R^N\times (0,T)\;  
\end{eqnarray*} 
since $u$ is a continuous solution of (\ref{slepcev}), $c^+ [u](x,t)\leq c[\psi^\alpha (u(\cdot,t)-u(x,t))]$ 
and $c^- [u](x,t)\geq c[\psi_\alpha (u(\cdot,t)-u(x,t))].$ 
The point for doing that is that the functions $c[\psi^\alpha (u(\cdot,t)-u(x,t))]$, $c[\psi_\alpha (u(\cdot,t)-u(x,t))]$ are now continuous in $x$ and $t$. 
 
Then we show that $\rho^+$, $\rho^-$ satisfy the same inequalities,  
the functions $c[\psi^\alpha (u(\cdot,t)-u(x,t))]$, $c[\psi_\alpha (u(\cdot,t)-u(x,t))]$  
being considered as fixed functions (in other words, we forget that they depend on $u$). 
In fact, we just provide the proof in details for $\rho^+$, the one for $\rho^-$ being analogous.  
Following the proof of \cite{bss93}, we set  
$$  
{u_\varepsilon} (x,t) := \frac{1}{2}\left(1 + \tanh \left( \varepsilon^{-1}(u(x,t) + \varepsilon^{1/2})\right)\right)\; . 
$$ 
Noticing that $u_\varepsilon = \phi_\varepsilon (u)$ for an increasing function $\phi_\varepsilon$, we have that 
the function ${u_\varepsilon}$ still satisfies the two above inequalities. It is easy to see that  
$$  
\rho^+ = \mathop{\rm lim\,sup\!^*\,} {u_\varepsilon} \quad \hbox{and} \quad  (\rho^+)_*  
= \mathop{\rm lim\,inf_*\,} {u_\varepsilon} \; , 
$$ 
and the half-relaxed limits method indeed shows that 
\begin{eqnarray*} 
(\rho^+)^*_t &\leq & c[\psi^\alpha (u(\cdot,t)-u(x,t))] |D(\rho^+)^*| \quad \mbox{in} 
\ \R^N\times (0,T)\; , \\ 
((\rho^+)_*)_t & \geq & c[\psi_\alpha (u(\cdot,t)-u(x,t))] |D(\rho^+)_*| \quad \mbox{in} 
\ \R^N\times (0,T)\; . 
\end{eqnarray*} 
 
The next step consists in remarking that the viscosity sub and 
supersolutions inequalities for $\rho^+$ are obviously satisfied in 
the complementary of $\partial \{u \geq 0\}$ since $\rho^+$ is locally 
constant there and therefore it is a classical solution of the problem. 
The only nontrivial viscosity sub and 
supersolutions inequalities we have to check are at points $(x,t)\in \partial \{u \geq 0\},$ 
i.e. such that $u(x,t)=0$ since $u$ is continuous. For such points, as 
$\alpha\to 0,$ 
\begin{eqnarray*} 
c[\psi^\alpha (u(\cdot,t)-u(x,t))]\to c[\rho^+](x,t)=  c[(\rho^+)^*](x,t) 
\end{eqnarray*} 
since $\rho^+$ is upper-semicontinuous, and 
\begin{eqnarray*} 
c[\psi_\alpha (u(\cdot,t)-u(x,t))]\to c[\rho^-](x,t). 
\end{eqnarray*} 
The stability result for equations with a $L^1$-dependence in time yields the inequalities 
\begin{eqnarray} 
\nonumber (\rho^+)^*_t &\leq& c[(\rho^+)^*]|D(\rho^+)^*| \quad \mbox{in} 
\ \R^N\times (0,T)\; , \\ 
\label{supsol1} ((\rho^+)_*)_t &\geq& c[\rho^-]|D(\rho^+)_*| \quad \mbox{in} 
\ \R^N\times (0,T)\; . 
\end{eqnarray} 
 
The second inequality is weaker than the one we claim : to obtain $c[(\rho^+)_*]$ 
instead of $c[\rho^-],$ we have to play with the different level-sets of $u$~: 
for $\beta >0$ small, we set $\rho^+_\beta = \1_{\{u \geq - \beta\}}$. Since 
$u$ is a solution of  (\ref{slepcev}) and $\psi_{\alpha, \beta}:= \psi_\alpha (\cdot +\beta)$ is nondecreasing, 
then $\psi_{\alpha, \beta} (u)$ is a (continuous) supersolution of 
\begin{eqnarray*} 
(\psi_{\alpha, \beta} (u))_t \geq c^- [\psi_{\alpha, \beta} (u)]| D\psi_{\alpha, \beta} (u)| 
 \quad \mbox{in} \ \R^N\times (0,T)\; . 
\end{eqnarray*} 
By stability we get, as $\alpha \to 0,$ 
\begin{eqnarray*} 
((\rho^+_\beta)_*)_t \geq c^- [(\rho^+_\beta)_*] |D(\rho^+_\beta)_*| 
\quad \mbox{in} \ \R^N\times (0,T)\; . 
\end{eqnarray*} 
But, for $(x,t)\in \partial \{u\geq -\beta\},$ 
\begin{equation} \label{ineq936} 
c^- [(\rho^+_\beta)_*](x,t) = c[\1_{ \{ (\rho^+_\beta)_*(\cdot ,t) > 0\} } ](x,t) 
= c[(\rho^+_\beta)_* ](x,t)\geq c[\rho^+](x,t). 
\end{equation} 
It follows that $(\rho^+_\beta)_*$ is a supersolution of the Eikonal Equation with 
$c[\rho^+](x,t)$ (as before, the only nontrivial inequalities we have to check 
are on $\partial \{u\geq -\beta\}$ and they are true because of (\ref{ineq936}). 
Letting $\beta$ tends to $0$ and using that $(\rho^+)_* =  \mathop{\rm lim\,inf_*\,} (\rho^+_\beta)_*,$ 
we obtain the expected inequality (even something better since (\ref{supsol1}) holds actually with 
$c[\rho^+](x,t)$). In particular, we get that $\rho^+$ is a solution of 
\begin{eqnarray} \label{eq::rajoutee} 
\left\{ 
\begin{array}{ll} 
\rho_t = c[\rho^+]|D\rho| & {\rm in  } \ \R^N\times [0,T)\; , \\[2mm] 
\rho(x,0) = \1_{\{u_0  \geq 0\}} & {\rm in  } \  \R^N\; . 
\end{array} 
\right. 
\end{eqnarray} 
And the proof of the claim is complete. \\ 
 
\noindent{\it 2. The functions $v^\pm$ are weak solutions of 
  (\ref{dislocation}).} 
Let us start with the ``$+$'' case. We 
first remark that the existence and uniqueness of $v^+$ follows from the 
standard theory for equations with a $L^1$-dependence in time 
(see Appendix A). 
 
To prove that $v^+$ is a weak solution of (\ref{dislocation}), it remains to prove 
that (\ref{FormeFaible2}) holds. It is sufficient to show that 
\begin{equation} \label{suite-incl} 
\{v^+ (\cdot, t) > 0\} \subset \{u(\cdot, t) \geq 0\} \subset \{v^+ (\cdot,t) \geq 0\} \; .  
\end{equation} 
We use again the functions $\psi_\alpha, \psi^\alpha$ introduced above. We remark that 
$$ \psi_\alpha (u_0) \leq \rho^+ (x,0) \leq \psi^\alpha (u_0)\quad \mbox{in  } \R^N\;.$$  
Moreover $v^+$ and $\rho^+$ are solutions of the same equation, namely 
(\ref{eq::rajoutee}) with $c[\rho^+]$ which is considered as a fixed function, and so are   
$\psi_\alpha (v^+)$ and $\psi^\alpha (v^+)$  
because the equation is geometric. Therefore, a standard comparison result implies 
$$  
\psi_\alpha (v^+) \leq (\rho^+)_* \leq  \rho^+ = (\rho^+)^* \leq \psi^\alpha (v^+)\quad  
\mbox{in  }  \R^N\times [0,T)\; .  
$$  
And letting $\alpha$ tends to $0$, these inequalities imply (\ref{suite-incl}). 
 
We can prove the symmetric result for $v^-$, the only difference is that inclusion (\ref{suite-incl}) has to be replaced by 
\begin{equation} \label{suite-incl-} 
\{v^- (\cdot, t) > 0\} \subset \{u(\cdot, t) > 0\} \subset \{v^- (\cdot,t) \geq 0\} \; .  
\end{equation} 
 
\noindent{\it 3. Claim: If $v$ is a weak solution of (\ref{dislocation}), then 
 $\1_{\{v (\cdot, t)\geq 0\}}$ is a $L^1$-subsolution of (\ref{eik111}).} 
From Proposition \ref{Trestresweak} and since $c_0\geq 0,$ $v$ satisfies in the $L^1$-sense, 
\begin{eqnarray} \label{sur-v1} 
v_t \leq c[\1_{\{v (\cdot, t) \geq 0\}}]|Dv| \quad \mbox{in  }  \R^N\times [0,T)\;. 
\end{eqnarray} 
By similar arguments as we used above, the function $\1_{\{v (\cdot, t)\geq 0\}}$ satisfies the same inequality 
which gives the result. \\ 
 
\noindent{\it 4. The function $\rho^+$ is the maximal $L^1$-subsolution of  (\ref{eik111}).} 
Let $w$ be a $L^1$ (upper-semicontinuous) subsolution of (\ref{eik111}).  
First we have $ w \leq 1$ in $\R^N\times [0,T)$ by comparison  
with the constant supersolution $1$ for the equation  
with $c[w]$ fixed. By considering $\max (w, 0)$ we may assume that  
$0 \leq w \leq 1$ in $\R^N\times [0,T)$.  
By similar arguments as we already used in Step 1, we can show that $\1_{\{w (\cdot, t) > 0\}}$ is also  
 a $L^1$-subsolution of (\ref{eik111}); thus we can assume that $w$ is a characteristic function. 
 
Then we remark that $w$ is also a subsolution of (\ref{slepcev}) : indeed, 
again, the only nontrivial viscosity inequalities are on the boundary of the 
set $\{w = 1\}$ and if $(x,t)$ is such a point we have $w(x,t)=1$ because $w$ is upper-semicontinuous  
and $w = \1_{\{w (\cdot, t) \geq w(x,t)\}}$.  
Since $u$ is a solution of the geometric equation (\ref{slepcev}), 
$\psi^\alpha(u)$ is still a solution which satisfies $\psi^\alpha(u)(x,0)\geq 
\1_{\{u_0 \geq 0\}}\geq w(x,0)$ in $\R^N.$ By Theorem \ref{comp-slep} we obtain 
$$  
w  \leq \psi^\alpha (u)\quad \mbox{in  }  \R^N\times [0,T)\; . 
$$ 
And letting $\alpha$ tend to $0$ provides $w\leq \rho^+$ which proves that  
$\rho^+$ is the maximal subsolution of (\ref{eik111}).\\ 
 
\noindent{\it 5. The function $v^+$ is the maximal weak solution of 
  (\ref{dislocation}).} 
Let $v$ be a weak solution of (\ref{dislocation}). 
From Steps 3 and 4, we get $\1_{\{v (\cdot, t) \geq 0\}} \leq \rho^+ (\cdot, t)$ in $\R^N\times [0,T)$ and 
(\ref{sur-v1}) implies 
$$  
v_t \leq c[\rho^+]|Dv| \quad \mbox{in  }  \R^N\times [0,T)\;. 
$$ 
Therefore $v$ is a subsolution of (\ref{nonlocaldisc})  
and by standard comparison result, this leads to $v \leq v^+$ in $\R^N\times [0,T),$ 
which proves the result. \\ 
 
\noindent{\it 6. We have $\{v^+ (\cdot,t) \geq 0\}= \{u(\cdot,t) \geq 0\}$ and 
$\{v^- (\cdot,t) > 0\}= \{u(\cdot,t) > 0\}.$} 
From Step 5, we get $\1_{\{v^+ (\cdot, t) \geq 0\}} \leq \rho^+ (\cdot, t) =\1_{\{u (\cdot, t)\geq 0\}}.$ The 
conclusion follows for $v^+$ using (\ref{suite-incl}). The inclusion for $v^-$ uses symmetric arguments.\\ 
 
\noindent{\it 7. Uniqueness when $\{u(\cdot ,t)=0\}$ has Lebesegue measure 0.} If 
$\mathcal{L}^N (\{u(\cdot ,t)=0\})=0,$ then $c[\rho^+]=c[\rho^-].$ Hence 
$v^+=v^-$ is the unique weak solution of (\ref{dislocation}) and it is obviously a classical one. \hfill $\Box$ 
 
\section*{Appendix A: a stability result for Eikonal Equations with 
  $L^1$-dependence in time}   
\label{appendice} 
 
The aim of this appendix is to provide a self-contained presentation of a 
stability result for viscosity 
solutions of Eikonal Equations with $L^1$-dependences in time which handles 
the case of weak convergence of the equations, instead of the classical 
strong $L^1$ convergence. This stability result is a particular case of a 
general stability result proved by Barles in \cite{barles06}. 
 
For $T>0$, we are interested in solutions of the following equation 
\begin{equation}\label{eq:cbar} 
\left\{\begin{array}{l} 
\displaystyle{\frac{\partial v}{\partial t}  = {\bar{c}}(x,t)|D v| \quad \mbox{dans} 
\quad  \R^N\times (0,T)}\\ 
\\ 
v(\cdot,0)=u_0\quad \mbox{dans  }\R^N\; , 
\end{array}\right. 
\end{equation} 
where the velocity ${\bar{c}} : \R^N\times (0,T) \to \R$ is defined for 
almost every $t\in (0,T)$. We 
also assume that ${\bar{c}}$ satisfies\\ 
\noindent {\bf (H3)} The function ${\bar{c}}$ is continuous with respect to $x\in \R^N$ 
and measurable in $t.$  
For all $x,y\in\R^N$ and almost all $t\in [0,T],$  
$$ 
|{\bar{c}}(x,t)|\le M \ \ \ {\rm and} \ \ \ 
|{\bar{c}}(x,t)-{\bar{c}}(y,t)|\le L|x-y|. 
$$ 
Let us underline that we do not assume any continuity in time of ${\bar{c}}$. 
We recall the following (under assumption {\bf (H0)}) 
\begin{Definition}{\bf ($L^1$-viscosity solutions)}\\ 
An upper-semicontinuous (respectively lower-semicontinuous) function $v$  on $\R^N\times [0,T]$  
is a $L^1$-viscosity subsolution (respectively supersolution) of 
(\ref{eq:cbar}), if  
$$ 
v(0,\cdot) \le u_0 \quad (\mbox{respectively}\quad v(0,\cdot) \ge u_0), 
$$ 
and if for every $(x_0,t_0)\in \R^N\times [0,T]$, $b\in L^1(0,T)$, $\varphi \in 
C^\infty(\R^N\times (0,T))$ and continuous function $G: \R^N\times (0,T) \times \R^N\to \R$ 
such that \\ 
(i) the function  
$$ 
(x,t)\longmapsto v(x,t)-\int_0^t b(s)ds -\varphi(x,t) 
$$  
has a local maximum 
(respectively minimum) at $(x_0,t_0)$ over $\R^N\times (0,T)$ and such that \\ 
(ii) for almost every $t\in (0,T)$ in some neighborhood of $t_0$ and for 
every $(x,p)$ in some neighborhood of $(x_0,p_0)$ with $p_0=\nabla\varphi(x_0,t_0)$, 
we have 
$${\bar{c}}(x,t)|p| -b(t)  \le G(x,t,p) \quad (\mbox{respectively}\quad 
{\bar{c}}(x,t)|p| -b(t)  \ge G(x,t,p))$$ 
then 
$$ 
\frac{\partial \varphi}{\partial t}(x_0,t_0) \le G(x_0,t_0,p_0)\quad 
(\mbox{respectively}\quad  \frac{\partial \varphi}{\partial t}(x_0,t_0) \ge 
G(x_0,t_0,p_0)). 
$$ 
Finally we say that a locally bounded function $v$ defined on $\R^N\times [0,T]$ 
is a $L^1$-viscosity solution of (\ref{eq:cbar}), if its upper-semicontinuous 
(respectively lower-semicontinuous) envelope is a 
$L^1$-viscosity subsolution (respectively supersolution). 
\end{Definition} 
Let us recall that viscosity solutions in the $L^1$-sense were introduced 
in Ishii's paper \cite{ishii85}.  
We refer to D. Nunziante \cite{nunziante90, nunziante92} and M.  
Bourgoing \cite{bourgoing04a,bourgoing04b} for a complete presentation of the theory. 
 
Then we have the following result 
 
\begin{Theorem}{\bf (Existence and uniqueness)}\\ 
For any $T>0$, under assumptions {\bf (H0)} and {\bf (H3)}, there exists a 
unique $L^1$-viscosity solution to (\ref{eq:cbar}). 
\end{Theorem} 
 
Finally, let us consider the solutions $v^\varepsilon$ to the following 
equation 
\begin{equation}\label{eq:cbarbis} 
\left\{\begin{array}{l} 
\displaystyle{\frac{\partial v^\varepsilon}{\partial t}  = {\bar{c}}^\varepsilon(x,t)|D v^\varepsilon| \quad \mbox{in} 
\  \R^N\times (0,T)},\\[2mm] 
v^\varepsilon(\cdot,0)=u_0\quad \mbox{in}\ \R^N\;. 
\end{array}\right. 
\end{equation} 
Then we have the following 
\begin{Theorem}\label{L1stability}{\bf ($L^1$-stability, \cite{barles06})}\\ 
Under assumption {\bf (H0)}, let us assume that the velocity 
${\bar{c}}^\varepsilon$ satisfies {\bf (H3)} (with some constants $M,L$ 
independent of $\varepsilon$). Let us consider the $L^1$-viscosity solution $v^\varepsilon$ to (\ref{eq:cbarbis}). Assume that 
$v^\varepsilon$ converges locally uniformly to a function $v$ and, for all $x\in \R^N,$ 
\begin{eqnarray} \label{cv-unif-t} 
\hspace*{0.5cm} \int_0^t {\bar{c}}^\varepsilon (x,s)ds \to \int_0^t {\bar{c}} (x,s)ds  
\ \ \ { locally \ uniformly \ in } \ (0,T). 
\end{eqnarray} 
Then $v$ is a $L^1$-viscosity solution of (\ref{eq:cbar}). 
\end{Theorem} 
 
\begin{Remark}\rm \label{rem-cvunif} 
Theorem \ref{L1stability} is stated as in \cite{barles06} but  
note that, under {\bf (H3)}, assumption (\ref{cv-unif-t}) is automatically 
satisfied as soon as the convergence is merely pointwise. Indeed, since 
$$ 
 \left| \int_0^t {\bar{c}}^\varepsilon (x,s)ds\right|\leq MT 
\  \ {\rm and} \ \  
\left| \int_0^t {\bar{c}}^\varepsilon (x,s)ds - \int_0^{t'} {\bar{c}}^\varepsilon (x,s)ds\right|  
\leq M |t-t'|, 
$$ 
from Ascoli's Theorem, the convergence is uniform. 
\end{Remark} 
 
\section*{Appendix B: Interior ball regularization (proof of Lemma \ref{lem:bornebouleint})}  \label{appendice-boul-int} 
 
The proof of this result can be adaptated from those of Cannarsa and 
Frankowska \cite{cf06} or \cite{acm05} (see also \cite{cc06} for related 
perimeter estimate for general equations). For the sake of completeness, we give a proof 
close to the one of \cite{cf06} (this latter holds for much more general, but time-independent, dynamics). 
The unique (and small) contribution of this part amounts to explain how this 
proof can be simplified in the particular case of dynamics of the  
form (\ref{EqDiff}) and to point out that the time-dependence is not an issue 
for the results of \cite{cf06} to hold.  
 
We first prove that the reachable set for controlled dynamics of the form  
\begin{equation}\label{EqDiff}  
\dot{x}(t)=c (x(t),t) u(t),\quad u \in L^\infty ([0,T],\overline{B}(0,1)) 
\end{equation}  
enjoys the interior ball property for positive time. We assume that 
$c:[0,T]\times \R^N\to \R$ satisfies, for any $x,y\in \R^N$ and $t\in [0,T],$  
\begin{equation*}\label{Hypf}  
\left\{\begin{array}{ll}  
(i) & c \;\mbox{\rm is Borel measurable,} \\ 
&  \mbox{differentiable  with respect to the space variable  for a.e. time,}\\  
(ii) & |c(x,t)-c(y,t)|\leq L_1 |x-y|,\\  
(iii)& |D_xc(x,t)-D_xc(y,t)|\leq N_1 |x-y|, \\  
(iv) & M_1\ge c(x,t)\geq \delta >0,  
\end{array}\right.  
\end{equation*}  
where $L_1,N_1\geq 0$ et $M_1,\delta>0$ are given constants. 
Let $K_0\subset \R^N$ be the initial set. 
We define the reachable set ${\mathcal R}(t)$ from $K_0$ for (\ref{EqDiff}) at 
time $t$ by:  
\begin{eqnarray*} 
{\mathcal R}(t) =\left\{ x(t) \;, \; x(\cdot) \; \mbox{\rm solution to (\ref{EqDiff}) with $x(0)\in K_0$} \right\}\;.  
\end{eqnarray*} 
It is known that ${\mathcal R}(t)$ is a closed subset of $\R^N$. Let $y_0$ be an extremal solution on the time interval $[0,T]$, i.e., a solution of (\ref{EqDiff}) such that  
$$  
y_0(0)\in K_0 \quad {\rm and }\quad y_0(T)\in \partial {\mathcal R}(T).  
$$  
From the Pontryagin Maximum Principle for extremal trajectories (see for 
instance \cite{clarke83}), 
there is some adjoint function $p_0:[0,T]\to \R^N\backslash\{0\}$ such that $(y_0,p_0)$ is a solution to:  
\begin{equation*}\label{SH}  
\left\{\begin{array}{l}  
\displaystyle \dot{y}_0(t)=c(y_0(t),t)\frac{p_0(t)}{|p_0(t)|},\\  
-\dot{p}_0(t)=D_xc(y_0(t),t)|p_0(t)|. 
\end{array}\right.  
\end{equation*}  
Since the system is positively homogeneous with respect to $p$, we can assume, without loss of generality, that  
$|p_0(T)|=1$ and we set $\theta_0:=p_0(T)$. 
 
Let $P$ be the matrix valued solution to  
\begin{equation*}  
\left\{\begin{array}{l}  
\displaystyle \dot{P}(t)=\frac{p_0(t)}{|p_0(t)|}[D_xc(y_0(t),t)]^* P(t),\\  
P(T)=Id. 
\end{array}\right.  
\end{equation*}  
A straightforward computation shows that $P^*(t)p_0(t)=\theta_0$ for any $t\in [0,T]$.  
 
Let us fix some parameter $\gamma>0$ to be chosen later, $\theta\in B(0,1)$ and 
let us set, for all $t\in [0,T],$ 
$$  
y_\theta(t)=y_0(t)-\gamma t P(t)(\theta_0-\theta).  
$$  
Our aim is to show that $y_\theta$ is a solution to (\ref{EqDiff}). Indeed  
we have  
\begin{eqnarray*} 
|\dot{y}_\theta|^2 & = & \left|c(y_0,t)\frac{p_0}{|p_0|} -  
\gamma P(\theta_0-\theta) - \gamma t \frac{p_0}{|p_0|} D_xc^* P(\theta_0-\theta) \right|^2 \\  
&= & |c(y_\theta,t)|^2 -2\gamma c(y_0,t) \left\langle \frac{p_0}{|p_0|} , P(\theta_0-\theta)\right\rangle\\  
&& + |c(y_0,t)|^2- |c(y_\theta,t)|^2 -2 \gamma t c(y_0,t)  \left\langle  \frac{p_0}{|p_0|},  \frac{p_0}{|p_0|}D_xc^*P(\theta_0-\theta) \right\rangle \\  
&& +\gamma^2 \left|P(\theta_0-\theta) + t \frac{p_0}{|p_0|} D_xc^* P(\theta_0-\theta) \right|^2\\  
&\leq & |c(y_\theta,t)|^2 -2\gamma \frac{c(y_0,t)}{|p_0|} 
\langle\theta_0,\theta_0-\theta\rangle \ \ \ ({\rm because} \ P^*p_0=\theta_0)\\  
&& + c^2(y_0,t)- c^2(y_\theta,t) - \langle D_x(c^2)(y_0,t), y_0-y_\theta \rangle +\gamma^2 M |\theta_0-\theta|^2\\  
&\leq & |c(y_\theta,t)|^2  
-\gamma \frac{\delta}{|p_0|} |\theta_0-\theta|^2\; +\; M_1'|y_0-y_\theta|^2+ \gamma^2M |\theta_0-\theta|^2\\  
&\leq & |c(y_\theta,t)|^2 -\gamma \frac{\delta}{|p_0|} 
 |\theta_0-\theta|^2\; +\; \gamma^2M' |\theta_0-\theta|^2,  
\end{eqnarray*} 
with $M_1'=L_1^2+M_1N_1$, and where $M$ and $M'$ only depend on $T, L_1, N_1,M_1$ because  
$|p_0(t)|$ is bounded from below by a constant depending only on $T, L_1$.  
Hence, for $\gamma$ sufficiently small, $y_\theta$ is a solution of (\ref{EqDiff}) starting from $y_0(0)\in K_0$ and therefore  
$y_\theta(T)\in {\mathcal R}(T)$. 
 
Finally ${\mathcal R}(T)$ contains all the $y_\theta(T)$ for $\theta\in 
B(0,1)$, i.e., the ball centered at $y_0(T)-\gamma T \theta_0$ and of 
radius $\gamma T$ (since $P(T)=Id$).  
 
We apply the previous result with $c =c_1$ and $K_0= \left\{v(\cdot,0)\ge 0\right\}.$  
Then  $\left\{v(\cdot,t)\ge 0\right\}= {\mathcal R}(t)$ for all $t>0.$ 
 
We end with a remark: in the statement of Lemma \ref{lem:bornebouleint}, 
$c_1$ is assumed to be continuous in time. As we have seen, it is not necessary, 
$c_1$ can be merely measurable in time up to consider the $L^1$-solution $v$ 
of (\ref{eq:2}) as recalled in Appendix A.  
\hfill $\Box$

\bigskip 
 
\noindent {\bf Aknowledgement }\\ 
This work was supported by the two contracts ACI JC 1025 (2003-2005) and ACI JC 1041 (2002-2004).


\end{document}